\documentclass{amsart}
\usepackage{amssymb,amscd}
\usepackage[alphabetic]{amsrefs}
\usepackage[all]{xy}

\numberwithin{equation}{subsection}

\swapnumbers
\theoremstyle{plain}
\newtheorem{thm}[subsection]{Theorem}
\newtheorem{prop}[subsection]{Proposition}
\newtheorem{cor}[subsection]{Corollary}
\newtheorem*{thm*}{Theorem}
\newtheorem{lemma}[subsection]{Lemma}

\theoremstyle{definition}

\theoremstyle{remark}
\newtheorem{rem}[subsection]{Remark}
\newtheorem{rems}[subsection]{Remarks}

\newcommand{\Curve}{\mathcal{C}}
\renewcommand{\O}{\mathcal{O}}
\newcommand{\FF}{\mathcal{F}}

\newcommand{\F}{\mathbb{F}}
\newcommand{\Fp}{{\mathbb{F}_p}}

\newcommand{\Fq}{{\mathbb{F}_q}}

\newcommand{\Z}{\mathbb{Z}}
\newcommand{\Q}{\mathbb{Q}}

\newcommand{\C}{\mathbb{C}}
\newcommand{\A}{\mathbb{A}}
\renewcommand{\P}{\mathbb{P}}

\renewcommand{\d}{\mathfrak{d}}

\newcommand{\n}{\mathfrak{n}}

\newcommand{\sha}{{\hbox to 10pt{\rlap{\hskip2.8pt\vrule
height6pt\hskip1.6pt\vrule height6pt\hskip1.6pt
\vrule height6pt}\hskip1pt\vrule height0.8pt width 8pt\hskip1pt}}}

\newcommand{\into}{\hookrightarrow}

\newcommand{\tensor}{\otimes}

\newcommand{\nodiv}{\not|}
\def\nodiv{\mathrel{\mathchoice{\not|}{\not|}{\kern-.2em\not\kern.2em|}{\kern-.2em\not\kern.2em|}}}

\newcommand{\GL}{\mathrm{GL}}
\newcommand{\Sp}{\mathrm{Sp}}

\DeclareMathOperator{\im}{Im}

\DeclareMathOperator{\res}{Res}

\DeclareMathOperator{\cond}{Cond}
\DeclareMathOperator{\ord}{ord}
\DeclareMathOperator{\rk}{Rank}

\DeclareMathOperator{\aut}{Aut}

\DeclareMathOperator{\gal}{Gal}
\DeclareMathOperator{\spec}{Spec}
\DeclareMathOperator{\en}{End}

\renewcommand{\vector}[1]{\mathbf{#1}}
\renewcommand{\v}{\vector{v}}

\newcommand{\Fpbar}{{\overline{\mathbb{F}}_p}}
\newcommand{\Fqbar}{{\overline{\mathbb{F}}_q}}

\newcommand{\Fr}{{\mathbb{F}}_{r}}
\newcommand{\Qbar}{{\overline{\mathbb{Q}}}}

\newcommand{\ratto}{{\dashrightarrow}}
\newcommand{\Ql}{{\mathbb{Q}_\ell}}

\newcommand{\Qlbar}{{\overline{\mathbb{Q}}_\ell}}
\DeclareMathOperator{\ind}{Ind}
\DeclareMathOperator{\sgn}{Sign}
\DeclareMathOperator{\swan}{Swan}

\newcommand{\XX}{\mathcal{X}}
\newcommand{\YY}{\mathcal{Y}}
\newcommand{\ZZ}{{\mathcal{Z}}}

\newcommand{\sdp}{{\rtimes}}

\begin{document}
\title[Large ranks over function fields]{$L$-functions with large analytic
rank and \\abelian varieties with large algebraic rank\\over function fields} 
\author{Douglas Ulmer}
\address{Department of Mathematics \\ University of Arizona \\ Tucson,
AZ  85721}
\email{ulmer@math.arizona.edu}
\thanks{This paper is based upon work partially supported by the National
Science Foundation under Grant No. DMS 0400877}
\date{July 17, 2006}
\maketitle

\section{Introduction}
The goal of this paper is to explain how a simple but apparently new
fact of linear algebra together with the cohomological interpretation
of $L$-functions allows one to produce many examples of $L$-functions
over function fields vanishing to high order at the center point of
their functional equation.  Conjectures of Birch and Swinnerton-Dyer,
Bloch, and Beilinson relate the orders of vanishing of some of these
$L$-functions to Mordell-Weil groups and other groups of algebraic
cycles.  For certain abelian varieties of high analytic rank, we are
also able to prove the conjecture of Birch and Swinnerton-Dyer thus
establishing the existence of large Mordell-Weil groups in those
cases.

In the rest of this section we state the main results of the paper.

\begin{thm}\label{thm:av1}
  For every prime number $p$, every positive integer $g$, and every
  integer $R$, there exist absolutely simple, non-isotrivial abelian
  varieties 
  $J$ of dimension $g$ over $\Fp(t)$ such that
  $\ord_{s=1}L(J/\Fp(t),s)\ge R$.  For all $p$ and $g$ there are
  examples of such $J$ for which the conjecture of Birch and
  Swinnerton-Dyer holds and so the rank of the finitely generated
  abelian group $J(\Fp(t))$ is at least $R$.
\end{thm}

The precise meaning of non-isotrivial and a stronger property enjoyed
by $J$ are explained in Section~\ref{ss:non-isotrivial}.

Our examples are completely explicit: We produce hyperelliptic curves
whose Jacobians have the properties asserted in the theorem.  For
example, if $p$ does not divide $(2g+2)(2g+1)$ then the Jacobian of
the curve with affine equation $y^2=x^{2g+2}+x^{2g+1}+t^{p^n+1}$ over
$\Fp(t)$ is absolutely simple, non-isotrivial, and has Mordell-Weil
group of rank $\ge p^n/2n$ over $\Fp(t)$.  This curve and similar
examples for other pairs $(p,g)$ meet asymptotic upper bounds on ranks
explained in Section~\ref{s:RankBounds}.

We can also produce high analytic ranks for $L$-functions of
cohomology groups of higher degree:

\begin{thm}\label{thm:av2}
  For every prime number $p>2$, every odd positive integer $k$, and
  every integer $R$, there exist infinitely many integers $g$ such
  that there exist absolutely simple, non-isotrivial abelian varieties
  $J$ of dimension $g$ over $\Fp(t)$ with
$$\ord_{s=(k+1)/2}L(H^k(J)_{prim},s)\ge R.$$
\end{thm}

Here $L(H^k(J)_{prim},s)$ is the $L$-function associated to the
primitive part of the $k$-th $\ell$-adic cohomology group of $J$.  See
Section~\ref{ss:la} for details. A conjecture (or rather ``recurring
fantasy'') of Bloch \cite{BlochRF} predicts that the order of
vanishing appearing in the theorem is equal to the rank of a group
of homologically trivial cycles of codimension $k$ on $J$ modulo
rational equivalence.  Producing the predicted cycles, even in
specific examples, looks like an interesting but difficult problem.

We also obtain new results on elliptic curves.  In \cite{UlmerR} large
ranks were obtained by considering a specific elliptic curve over
various rational extensions of the base field.  The following result
shows that this is a very general phenomenon.

\begin{thm}\label{thm:ec1}
  Let $\Fq$ be the field with $q$ elements, $q$ a power of $p$, and
  let $E$ be any elliptic curve defined over $F=\Fq(v)$ such that the
  $j$-invariant of $E$ does not lie in $\Fq$.  Then there exists a
  power $r$ of $q$ such that for every integer $R$ there are
  extensions of $F$ of the form $K=\F_{r}(t)$ such that
  $\ord_{s=1}L(E/K,s)\ge R$.
\end{thm}

Regarding isotrivial elliptic curves, our method also gives a new
proof of a result of Tate and Shafarevitch:

\begin{thm}\label{thm:ec2}
  Let $E_0$ be a supersingular elliptic curve over $\Fp$ and let
  $E=E_0\times_{\spec\Fp}\spec\Fp(t)$.  Then for every integer $R$
  there exist quadratic twists $E'$ of $E$ over $\Fp(t)$ such that the
  rank of $E'(\Fp(t))$ is $\ge R$.
\end{thm}


\subsection{}
The key result of linear algebra and its implications for
$L$-functions already appeared in our previous work \cite{UlmerGNV} on
non-vanishing of $L$-functions.  (In that context, it was something of
a technicality, but here it returns in a more appealing guise.)  For
the convenience of the reader, we give a brief review of the linear
algebra from a somewhat different point of view and a more general
application to $L$-functions in Sections~\ref{s:la} through
\ref{s:L}. We then prove the results stated above in
Sections~\ref{s:av1} through \ref{s:ec2}.  In
Section~\ref{s:RankBounds} we discuss an upper bound on ranks in terms
of conductors and then note that the results of
Section~\ref{s:algranks} show that the main term of the bound is
sharp.

Our analytic rank results are all based on an understanding of the
behavior of $L$-functions in towers of function fields, the simplest
and most important example being the tower $\Fq(t^{1/d})$ where $d$
runs over integers prime to $p$, the characteristic of $\Fq$.  That
ranks of $L$-functions should often be unbounded in towers became
apparent while considering a question of Ellenberg on towers over
finite fields versus towers over number fields.  In a companion
\cite{UlmerR3} to this paper, we explain Ellenberg's question and
ultimately answer it in the negative by giving several examples of
abelian varieties which have ranks over $\Fq(t^{1/d})$ bounded
independently of $d$.

\subsection{}
It is a pleasure to thank Jordan Ellenberg for his stimulating
questions about ranks of elliptic curves as well as Brian Conrey, Bill
McCallum, Dinesh Thakur, and especially Bjorn Poonen for their help.

\section{Linear algebra}\label{s:la}
\begin{prop}\label{thm:la}
  Let $V$ be a finite-dimensional vector space with subspaces $W_i$
  indexed by $i\in\Z/a\Z$ such that $V=\oplus_{i\in\Z/a\Z}W_i$.  Let
  $\phi:V\to V$ be an invertible linear transformation such that
  $\phi(W_i)= W_{i+1}$ for all $i\in\Z/a\Z$.  Suppose that $V$ admits
  a non-degenerate, $\phi$-invariant bilinear form $\langle ,\rangle $
  which is either symmetric \textup{(}in which case we set
  $\epsilon=1$\textup{)} or skew-symmetric \textup{(}in which case
  $\epsilon=-1$\textup{)}.  Suppose that $a$ is even and
  $\langle,\rangle$ induces an isomorphism $W_{a/2}\cong W_0^*$
  \textup{(}the dual vector space of $W_0$\textup{)}.  Suppose also
  that $N=\dim W_0$ is odd.  Then the polynomial $1-\epsilon T^{a}$
  divides $\det(1-\phi T|V)$.
\end{prop}

The proof of Proposition~\ref{thm:la} is given in
Subsections~\ref{lemma:cyclic} through \ref{ss:la-end} below.

\begin{lemma}  
\label{lemma:cyclic}
Let $V$ be a finite-dimensional vector space with subspaces $W_i$
indexed by $i\in\Z/a\Z$ such that $V=\oplus_{i\in\Z/a\Z}W_i$.  Let
$\phi:V\to V$ be a linear transformation such that $\phi(W_i)\subset
W_{i+1}$ for all $i\in\Z/a\Z$.  Then
\begin{equation}\label{eqn:charpoly}
\det(1-\phi T|V)=\det(1-\phi^{a}T^a|W_0)
\end{equation}
\end{lemma}

\begin{proof}
  We argue by induction on the dimension of $W_0$.  If $W_0=\{0\}$
  then $\phi$ is nilpotent and both sides of \ref{eqn:charpoly} are 1.
  We may assume that the ground field is algebraically closed and so
  if $W_0\neq\{0\}$ then $\phi^a$ has an eigenvector $v\in W_0$ with
  eigenvalue $\lambda$.  Let $W_i'$ be the span of $\phi^iv$ and
  $V'=\oplus W_i'$.  A simple computation shows that
$$\det(1-\phi T|V')=1-\lambda T^a=\det(1-\phi^aT^a|W_0').$$
Since characteristic polynomials are multiplicative in short exact
sequences, we may replace $V$ with $V/V'$ and $W_i$ with $W_i/W_i'$
and finish by induction on the dimension of $W_0$.
\end{proof}

\begin{lemma}\label{lemma:eigenvalues}
  Under the hypotheses of Proposition~\ref{thm:la}, if $\lambda$ is an
  eigenvalue of $\phi^{a}$ on $W_0$ then so is $\lambda^{-1}$.
\end{lemma}

\begin{proof}
  First we note that since $\phi^{a}:W_0\to W_0$ factors as
\begin{equation*}
\xymatrix{W_0\ar[r]^{\phi^{a/2}}&W_{a/2}\ar[r]^{\phi^{a}}
&W_{a/2}\ar[r]^{\phi^{-a/2}}&W_0}
\end{equation*}
the eigenvalues of $\phi^{a}$ on $W_0$ are the same as the eigenvalues
of $\phi^{a}$ on $W_{a/2}$.  On the other hand, the pairing $\langle
,\rangle $ induces a duality between $W_0$ and $W_{a/2}$ for which
$\phi^{a}$ is orthogonal (i.e., for all $v\in W_0$, $w\in W_{a/2}$,
$\langle \phi^{a}v,w\rangle=\langle v,\phi^{-a}w\rangle$) and so the
eigenvalues of $\phi^{a}$ on $W_0$ are the inverses of the eigenvalues
of $\phi^{a}$ on $W_{a/2}$.
\end{proof}

\begin{lemma}\label{lemma:det}
  Under the hypotheses of Proposition~\ref{thm:la}, the determinant of
  $\phi^{a}:W_0\to W_0$ is $\epsilon^N$.
\end{lemma}

\begin{proof}
  The pairing $\langle ,\rangle $ induces a pairing on
  $W=\bigwedge^NW_0\oplus\bigwedge^NW_{a/2}$ which we again denote by
  $\langle ,\rangle $.  The sign of this pairing is $\epsilon^N$,
  i.e., $\langle v,w\rangle =\epsilon^N\langle w,v\rangle $ for all
  $v,w\in W$.  Let $h:W\to W$ be induced by $\bigwedge^N\phi^{a/2}$ and note
  that $h$ exchanges the subspaces $\bigwedge^NW_0$ and
  $\bigwedge^NW_{a/2}$.  Choose $v\in\bigwedge^NW_0$ and
  $w\in\bigwedge^NW_{a/2}$ such that $\langle v,w\rangle =1$.  Then
$$\det(\phi^{a}|W_0)=\langle h^2v,w\rangle =\langle hv,h^{-1}w\rangle 
=\langle w,v\rangle =\epsilon^N\langle v,w\rangle =\epsilon^N.$$
\end{proof}

\subsection{}\label{ss:la-end}
Proposition~\ref{thm:la} is an easy consequence of the lemmas.
Indeed, Lemmas~\ref{lemma:eigenvalues} and \ref{lemma:det} imply that
$\epsilon$ is an eigenvalue of $\phi^{a}$ on $W_0$, i.e., that
$1-\epsilon T$ divides $\det(1-\phi^{a}T|W_0)$ and then
Lemma~\ref{lemma:cyclic} implies that $1-\epsilon T^{a}$ divides
$\det(1-\phi T|V)$.  This completes the proof of
Proposition~\ref{thm:la}.

\begin{rems}\hfil\break
\begin{enumerate}
\item \vskip-12pt Under the hypotheses of Proposition~\ref{thm:la},
  $\epsilon\phi^a$ is the asymmetry (in the sense of
  \cite{Cortella-Tignol}) of the pairing $(w,w')=\langle
  w,\phi^{-a/2}w'\rangle$ on $W_0$.  This provides another way to see
  that $\epsilon\phi^a$ has determinant 1 and is conjugate to its
  inverse.
\item With hypotheses as in Proposition~\ref{thm:la} except with $N$
  even, we do not get any consequences for the eigenvalues of $\phi$
  except what is forced by Lemma~\ref{lemma:eigenvalues}.  See
  \cite{UlmerGNV}*{7.1.12} for a more precise version of this remark.
\end{enumerate}\end{rems}

On the other hand, combining Lemma~\ref{lemma:cyclic} with well-known
facts about orthogonal transformations yields the following variant,
whose proof will be left to the reader.

\begin{prop}\label{thm:la-var}
  Let $V$ be a finite-dimensional vector space with subspaces $W_i$
  indexed by $i\in\Z/a\Z$ such that $V=\oplus_{i\in\Z/a\Z}W_i$.  Let
  $\phi:V\to V$ be an invertible linear transformation such that
  $\phi(W_i)= W_{i+1}$ for all $i\in\Z/a\Z$.  Suppose that $V$ admits
  a $\phi$-invariant bilinear form $\langle ,\rangle $ such that
  $\langle,\rangle$ restricted to $W_0$ is non-degenerate and
  symmetric.  If $N=\dim W_0$ is odd and $\epsilon=\det(\phi^a|W_0)$,
  then $1-\epsilon T^{a}$ divides $\det(1-\phi T|V)$.  If $N=\dim W_0$
  is even and $\det(\phi^a|W_0)=-1$, then $1- T^{2a}$ divides
  $\det(1-\phi T|V)$.
\end{prop}

\section{Group theory}\label{s:gt}
We review some simple facts about the representation theory of an
extension of a finite abelian group by a cyclic group.  Fix an
algebraically closed field $k$ of characteristic zero.  In the
applications, $k$ will be $\Qlbar$.

\subsection{}
Let $H$ be a finite abelian group and let $\phi:H\to H$ be an
automorphism of $H$.  Let $C$ be the cyclic subgroup of $\aut(H)$
generated by $\phi$ and let $b$ denote the order of $C$.  We form the
semidirect product $H^+=H\sdp C$; explicitly, $H^+$ is the set of
pairs $(h,\phi^i)$ with $h\in H$ and $i\in\Z/b\Z$ with multiplication
$(h,\phi^i)(h',\phi^j)=(h\phi^i(h'),\phi^{i+j})$.  For $a$ an integer,
let $H^+_a$ be the subgroup of $H^+$ generated by $H$ and $\phi^a$; it
has index $\gcd(a,b)$ in $H^+$.

\subsection{}
Let $\hat H$ denote the group of $k^\times$-valued characters of $H$.
There is a natural action of $C$ on $\hat H$: if $\chi\in\hat H$ and
$h\in H$, then $\chi^\phi(h)$ is defined to be $\chi(\phi(h))$.  Given
$\chi\in\hat H$, let $a=a_\phi$ be the smallest positive integer such
that $\chi^{\phi^a}=\chi$.  Choose a $(b/a)$-th root of unity
$\zeta\in k$.  We extend $\chi$ to a character $\tilde \chi$ of
$H^+_a$ by setting $\tilde\chi(\phi^a)=\zeta$.  It is not hard to
check that the induced representation $\ind^{H^+}_{H^+_a}\tilde\chi$
is irreducible and up to isomorphism it only depends on $\zeta$ and
the orbit of the $C$ action on $\hat H$ containing $\chi$.  We denote
this orbit by $o$ and write $\sigma_{o,\zeta}$ for
$\ind^{H^+}_{H^+_a}\tilde\chi$.  Every irreducible
representation of $H^+$ is isomorphic to a $\sigma_{o,\zeta}$ for a
unique pair $(o,\zeta)$.  (This is a special case of the ``method of
little groups.''  See \cite{SerreLRFG}*{8.2} for details.)

It is easy to see that the dual of $\sigma_{o,\zeta}$ is
$\sigma_{-o,\zeta^{-1}}$ where $-o=\{\chi^{-1}|\chi\in o\}$.  In
particular, $\sigma_{o,\zeta}$ is self-dual if and only if $o=-o$ and
$\zeta\in\{\pm1\}$.  We write $\sigma_o$ for $\sigma_{o,1}$.

\subsection{}
Let $\Sigma=\ind^{H^+}_{C}\boldsymbol1$ where we write $\boldsymbol1$
for the trivial representation of $C$ with coefficients in $k$.  I
claim that
\begin{equation}\label{eq:sigma-decomp}
\Sigma\cong\bigoplus_{o\subset\hat H}\sigma_o
\end{equation}
where the sum on the right is over the
orbits of $C$ acting on $\hat H$.  Indeed, for each pair $(o,\zeta)$
choose $\chi\in o$ and extend it to $\tilde\chi$ as above.  Then using
standard notation for the inner product on the representation rings of
$H^+$ and $C$, we have
\begin{align*}
\langle\sigma_{o,\zeta},\Sigma\rangle_{H^+}
&=\langle\ind^{H^+}_{H^+_a}\tilde\chi,\ind^{H^+}_{C}\boldsymbol1\rangle_{H^+}\cr
&=\langle\res^{H^+}_{C}\ind^{H^+}_{H^+_a}\tilde\chi,\boldsymbol1\rangle_{C}
\end{align*}
which is the multiplicity of $1$ as eigenvalue of $\phi$ on
$\ind^{H^+}_{H^+_a}\tilde\chi$.  By Lemma~\ref{lemma:cyclic}, this is
the same as the multiplicity of 1 as an eigenvalue of
$\tilde\chi(\phi^a)$, namely 1 if $\zeta=1$ and 0 if $\zeta\neq1$.
Thus each $\sigma_{o}$ appears in $\Sigma$ exactly once and no
$\sigma_{o,\zeta}$ with $\zeta\neq1$ appears.  This establishes the
claim.

\subsection{}
Note that $(\chi^\phi)^{-1}=(\chi^{-1})^\phi$ so if $o\subset\hat H$
is an orbit of $C$ such that $o=-o$ (a ``self-dual orbit'') then the
involution $\chi\mapsto\chi^{-1}$ of $o$ is either trivial or has no
fixed points.  The first case happens exactly when $o$ consists
entirely of characters of order dividing 2 (in which case we say that
$o$ ``consists of order 2 characters'') and the second case happens
when all the characters in $o$ have (the same) order larger than 2 (in
which case we say $o$ ``consists of higher order characters'').

\subsection{}\label{ss:genl-app}
In the applications of these results, $F$ will be the function field
of a curve over a finite field $\Fq$, and $K$ will be a finite
extension of $F$ which is ``geometrically abelian,'' i.e., such that
the extension $\Fqbar K/\Fqbar F$ is abelian.  Then $H$ will be
$\gal(\Fqbar K/\Fqbar F)$ and $\phi$ will be the action of the
geometric ($q^{-1}$-power) Frobenius on $H$.  It is easy to see that
$b$ is then the degree of the algebraic closure of $\Fq$ in the Galois
closure $L$ of $K/F$ and we have the diagram of fields
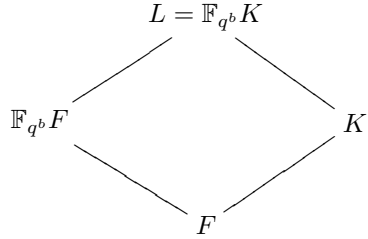
\begin{equation*}
\xymatrix{&L=\F_{q^b}K\ar@{-}[dl]\ar@{-}[dr]&\\
\F_{q^b}F\ar@{-}[dr]&&K\ar@{-}[dl]\\
&F&}
\end{equation*}
and the corresponding diagram of Galois groups:
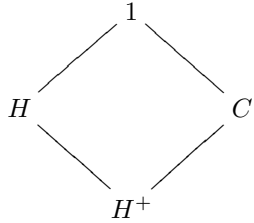
\begin{equation*}
\xymatrix{&1\ar@{-}[dl]\ar@{-}[dr]&\\
H\ar@{-}[dr]&&C\ar@{-}[dl]\\
&H^+&}
\end{equation*}

\subsection{}\label{ss:special-d}
Specializing further, the most interesting applications will be
in the case where $d$ is an integer prime to the characteristic of $F$
and $K=F(u^{1/d})$ for some $u\in F$ such that $[K:F]=d$.  In this
case, $H=\mu_d$ by Kummer theory, $\hat H=\Z/d\Z$, and the action of
$\phi$ is just multiplication by $q^{-1}$.  There are at most two orbits
$o$ consisting of characters of order 2, namely $o=\{0\}$ and, if $d$
is even, $o=\{d/2\}$.  On the other hand there is a plentiful supply
of self-dual orbits consisting of higher order characters.  Indeed, if
$d$ divides $q^n+1$ for some $n$, then $q^n\equiv-1\pmod{d}$ and so
every orbit is self-dual.  Since $q^{2n}\equiv1\pmod{d}$ each orbit
has cardinality at most $2n$ and so there are at least $(q^n-1)/2n$
self-dual orbits consisting of higher order characters.

\subsection{}
The results of this section can be extended, with some additional
complications, to the case where $H$ is an arbitrary finite group and
the extended results seem to have interesting applications to
arithmetic.  I hope to report on this elsewhere.

\section{Application to $L$-functions}\label{s:L}
We now apply the linear algebra result Proposition~\ref{thm:la} to
$L$-functions.  The discussion is a generalization of
\cite{UlmerGNV}*{3.2, 4.2, and 7.1}.

\subsection{}\label{ss:F}
Let $\Curve$ be a smooth, proper, geometrically irreducible curve over
the finite field $\Fq$ of characteristic $p$ and let $F=\Fq(\Curve)$
be its field of functions.  Choose an algebraic closure $F^{\rm alg}$
of $F$ and let $\overline{F}\subset F^{\rm alg}$ be the separable
closure of $F$.  Let $G_F=\gal(\overline{F}/F)$ be the absolute Galois
group of $F$.  For each place $v$ of $F$ we choose a decomposition
group $D_v\subset G_F$ and we let $I_v$ and $Fr_v$ be the corresponding
inertia group and geometric Frobenius class.  We write $\deg v$ for
the degree of $v$ and $q_v=q^{\deg v}$ for the cardinality of the
residue field at $v$.  For a finite extension $K$ of $F$, we denote
$\gal(\overline F/K)$ by $G_K$.

Fix a prime $\ell\neq p$ and let $\Qlbar$ be an algebraic closure
of $\Ql$, the field of $\ell$-adic numbers.  Fix also imbeddings
$\overline\Q\into\C$ and $\overline\Q\into\Qlbar$ and a compatible
isomorphism $\iota:\Qlbar\to\C$.  Whenever a square root of $q$ is
needed in $\Qlbar$, we take the one mapping to the positive square
root of $q$ in $\C$.  Having made this choice, we can define Tate
twists by half integers.

\subsection{}\label{ss:rho-hyps}
Fix a continuous representation
$\rho:G_F\to\GL_r(\Qlbar)$.
(As is well-known, $\rho$ factors through $\GL_r(E)$ for some finite
extension $E$ of $\Ql$.  See \cite{KS}*{9.0.7-9.0.8} for a 
discussion.) 
We assume that $\rho$ satisfies the following conditions:
\begin{enumerate}
\item $\rho$ is unramified outside a finite set of places, so that it
  factors through $\pi_1(U,\overline{\eta})$ for some non-empty open
  subscheme $j:U\into\Curve$.  (Here $\overline{\eta}$ is the
  geometric point of $\Curve$ defined by the fixed embedding $F\into
  F^{\rm alg}$.)
\item $\rho$ is $\iota$-pure of some integer weight $w$, i.e., for
  every place $v$ where $\rho$ is unramified, each eigenvalue $\alpha$
  of $\rho(Fr_v)$ satisfies $|\iota(\alpha)|=q_v^{w/2}$.
\item $\rho$ is self-dual of weight $w$ and sign
  $\sgn(\rho)\in\{\pm1\}$.  In other words, we assume that the space
  $\Qlbar^r$ on which $G_F$ acts via $\rho$ admits a non-degenerate,
  $G_F$-equivariant bilinear pairing $\langle ,\rangle $ with values in
  $\Qlbar(-w)$ and with $\langle v,v'\rangle
  =\sgn(\rho)\langle v',v\rangle $ for all $v,v'\in \Qlbar^r$.
\end{enumerate}

For each place $v$ of $F$ we write $\cond_v\rho$ for the exponent of 
the Artin conductor of $\rho$ at $v$.  (See~\cite{SerreLF}*{Chap.~VI}
for definitions.)  We let $\cond(\rho)=\sum_v({\cond_v\rho})[v]$ be the global
Artin conductor of $\rho$, viewed as an effective divisor on $\Curve$.

\subsection{}\label{ss:L-funs}
Attached to $\rho$ we have an $L$-function, defined formally by a
product over the places of $F$: 
\begin{align*}
L(\rho,F,T)&=\prod_v\det\left(1-\rho(Fr_v)T^{\deg v}
                  \left|(\Qlbar^r)^{\rho(I_v)}\right.\right)^{-1}
\end{align*}
and, for a complex variable $s$, we define $L(\rho,F,s)$ to be
$L(\rho,F,q^{-s})$.

Grothendieck's analysis of $L$-functions shows that $L(\rho,F,T)$ is a
rational function in $T$ and satisfies the functional equation
\begin{equation*}
L(\rho,F,T)=\left(q^{\frac{w+1}2}T\right)^NL(\rho,F,(q^{w+1}T)^{-1})
\end{equation*}
where
\begin{equation*}
N=(2g_\Curve-2)(\deg\rho)+\deg(\cond(\rho)).
\end{equation*}
(See, for example, \cite{MilneEC}, especially Section~VI.13.)

If $K$ is a finite extension of $F$ contained in $\overline{F}$, we
abbreviate $L(\rho|_{G_K},K,T)$ to $L(\rho,K,T)$.

\subsection{}\label{ss:geom-abelian-base-change}
Fix a finite extension $K$ of $F$ which is geometrically abelian in
the sense that $\Fqbar K/\Fqbar F$ is Galois.  We adopt the
definitions and notation of Subsection~\ref{ss:genl-app}, so that 
the Galois closure of $K/F$ is $L=\F_{q^b}K$, $H=\gal(\Fqbar K/\Fqbar F)$,
$C=\gal(\F_{q^b}K/K)\cong\gal(\F_{q^b}/\Fq)$ generated by the $q^{-1}$-power
Frobenius $\phi$, and $H^+=\gal(L/F)\cong H\sdp C$.

Continuing with the notations of Section~\ref{s:gt} we let
$\Sigma=\ind^{H^+}_{C}{\boldsymbol 1}$ so that
$\Sigma\cong\oplus_{o\subset \hat H}\sigma_o$ where the sum is over
orbits of $C$ on $\hat H$, the dual group of $H$.  We view $\Sigma$
and the $\sigma_o$ as representations of $G_F$ via the natural
surjection $G_F\to H^+$.

Now consider 
$L(\rho,K,T)=L(\res^{G_F}_{G_K}\rho,K,T)$.  
By standard properties of $L$-functions (e.g.,
\cite{DeligneConstants}*{3.8}) and basic representation theory,
\begin{align*}
L(\res^{G_F}_{G_K}\rho,K,T)
&=L(\ind^{G_F}_{G_K}
     \res^{G_F}_{G_K}\rho,F,T)\cr
&=L(\rho\tensor
        \ind^{G_F}_{G_K}\boldsymbol1,F,T)\cr
&=L(\rho\tensor\Sigma,F,T)\cr
&=\prod_{o\subset\hat H}L(\rho\tensor\sigma_o,F,T)
\end{align*}

Our basic result about $L$-functions says that for a suitable $K$,
many of the factors on the right hand side of the last equation vanish
at the center point of their functional equations:

\begin{thm}\label{thm:Lzeroes}
  Let $F$, $\rho$, and $K$ be as in \ref{ss:F}, \ref{ss:rho-hyps}, and
  \ref{ss:geom-abelian-base-change} respectively.  We keep the
  notations $H$, $\hat H$, and $C$ of
  \ref{ss:geom-abelian-base-change}.  Fix an orbit $o\subset\hat H$
  for the action of $C$ of cardinality $|o|$ which is self-dual
  \textup{(}$o=-o$\textup{)} and consists of characters of higher
  order \textup{(}$\chi\in
  o\implies\chi\neq\chi^{-1}$\textup{)}. Assume that for one
  \textup{(}and thus every\textup{)} $\chi\in o$ the degree of
  $\cond(\rho\tensor\chi)$ is odd.  Let $w$ be the weight of $\rho$
  and let $\epsilon=-\sgn(\rho)$.  Then
  $1-\epsilon\left(Tq^{\frac{w+1}2}\right)^{|o|}$ divides the
  numerator of $L(\rho\tensor\sigma_{o},F,T)$.
\end{thm}

\begin{proof}
  The theorem is a fairly straightforward application of the linear
  algebra result of Section~\ref{s:la} and the cohomological
  interpretation of $L$-functions.

  Let $j:U\into\Curve$ be a non-empty open subscheme over which both
  $\rho$ and $\sigma_{o}$ (and therefore also
  $\rho\tensor\sigma_{o}$) are unramified. These three
  representations give rise to lisse $\ell$-adic sheaves on $U$ and we
  let $\FF_\rho$, $\FF_{\sigma_{o}}$, and
  $\FF_{\rho\tensor\sigma_{o}}$ denote their direct images under $j$
  on $\Curve$.  (These are the ``middle extension'' sheaves on
  $\Curve$ attached to the representations.)

  Grothendieck's analysis of $L$-functions and our hypotheses on
  $\rho$ give a cohomological calculation of
  $L(\rho\tensor\sigma_{o},F,T)$:
\begin{equation*}
L(\rho\tensor\sigma_{o},F,T)
=\frac{\det\left(1-\phi T|H^1(\Curve\times\Fqbar,\FF_{\rho\tensor\sigma_{o}})\right)}{\det\left(1-\phi T|H^0(\Curve\times\Fqbar,\FF_{\rho\tensor\sigma_{o}})\right)\det\left(1-\phi T|H^2(\Curve\times\Fqbar,\FF_{\rho\tensor\sigma_{o}})\right)}
\end{equation*}
where $\phi$ is the geometric Frobenius in $\gal(\Fqbar/\Fq)$.  By
Deligne's theorem on weights, there is no cancellation in this
expression and so the numerator of the $L$-function is precisely
\begin{equation*}
\det\left(1-\phi
  T|H^1(\Curve\times\Fqbar,\FF_{\rho\tensor\sigma_{o}})\right).
\end{equation*}

The theorem is invariant under twisting and so we may replace $\rho$
with $\rho\tensor \Qlbar(\frac{w+1}2)$ and assume that $\rho$ is self-dual
of weight $-1$ and sign $\sgn(\rho)$.  Since $-o=o$, $\sigma_{o}$ is
self-dual with sign +1 and so $\rho\tensor\sigma_{o}$ is self-dual
with sign $\sgn(\rho)$.  Poincar\'e duality implies that
$H^1(\Curve\times\Fqbar,\FF_{\rho\tensor\sigma_{o}})$ is self-dual
of weight 0 and sign $\epsilon=-\sgn(\rho)$ as a representation of
$\gal(\Fqbar/\Fq)$.

On the other hand, $\sigma_{o}$ factors as a representation of
$\gal(\overline{F}/\Fqbar F)$ into lines and so on
$\overline{\Curve}=\Curve\times\Fqbar$
\begin{equation*}
\FF_{\rho\tensor\sigma_{o}}\cong\bigoplus_{\chi\in o}\FF_{\rho\tensor\chi}
\end{equation*}
where $\FF_{\rho\tensor\chi}$ is the middle extension sheaf attached
to $\rho\tensor\chi$.  Thus we have a factorization
\begin{equation*}
H^1(\Curve\times\Fqbar,\FF_{\rho\tensor\sigma_{o}})
=\bigoplus_{\chi\in o}
H^1(\Curve\times\Fqbar,\FF_{\rho\tensor\chi}).
\end{equation*}
Under the Poincar\'e duality pairing on
$H^1(\Curve\times\Fqbar,\FF_{\rho\tensor\sigma_{o}})$, the subspaces
$H^1(\Curve\times\Fqbar,\FF_{\rho\tensor\chi})$ and
$H^1(\Curve\times\Fqbar,\FF_{\rho\tensor\chi^{-1}})$ are dual to one another.
Moreover, $\phi$ preserves the pairing and sends
$H^1(\Curve\times\Fqbar,\FF_{\rho\tensor\chi})$ to
$H^1(\Curve\times\Fqbar,\FF_{\rho\tensor\chi^{\phi}})$.  

For any middle extension sheaf $\FF_\tau$ on $\Curve$ associated to an
$\ell$-adic representation $\tau$ of $G_F$ 
satisfying the first hypothesis of Subsection~\ref{ss:rho-hyps}, we have 
\begin{equation*}
\dim H^0(\Curve\times\Fqbar,\FF_{\tau})=\dim
H^2(\Curve\times\Fqbar,\FF_{\tau})
\end{equation*}
both of these being the multiplicity with which the trivial
representation appears in $\tau$ restricted to $\gal(\overline
F/\Fqbar F)$.  It follows that the dimension of
$H^1(\Curve\times\Fqbar,\FF_{\tau})$ has the same parity as the Euler
characteristic $(2-2g_\Curve)\deg(\tau)-\deg\cond(\tau)$ and this has
the same parity as the degree of $\cond(\tau)$.  Therefore, 
our hypotheses imply that the dimension of
$H^1(\Curve\times\Fqbar,\FF_{\rho\tensor\chi})$ is odd. 

Theorem~\ref{thm:Lzeroes} now follows easily from
Proposition~\ref{thm:la}.  Indeed, fix $\chi\in o$ and set
$V=H^1(\Curve\times\Fqbar,\FF_{\rho\tensor\sigma_{o}})$ and
$W_i=H^1(\Curve\times\Fqbar,\FF_{\rho\tensor\chi^{\phi^i}})$.  The
geometric Frobenius $\phi$ permutes the $W_i$ cyclically.  Since $o$
is self-dual consisting of higher order characters, $a=|o|$ is even.
Poincar\'e duality give a non-degenerate pairing on $V$ which induces
a duality between $W_0$ and $W_{a/2}$ and the dimension of
$W_0$ is odd.  Thus the hypotheses of Proposition~\ref{thm:la} are
satisfied and so $1-\epsilon T^{a}$ divides the numerator of the
twisted $L$-function and $1-\epsilon (q^{\frac{w+1}2}T)^{a}$ divides
the numerator of the untwisted $L$-function
$L(\rho\tensor\sigma_{o},F,T)$.
\end{proof}

\begin{rem}
One can formulate a variant of Theorem~\ref{thm:Lzeroes} with 
Proposition~\ref{thm:la-var} playing the role of
Proposition~\ref{thm:la}.  This variant does not seem to lead to
unbounded ranks and so we omit it.
\end{rem}

We now give a context in which Theorem~\ref{thm:Lzeroes} can be
applied to deduce unbounded ranks in towers.  For the definition of
the Swan conductor of a representation, we refer to
\cite{MilneEC}*{p.~188}. 

\begin{thm}\label{thm:towers}
  Let $F=\Fq(u)$ where $q$ is a power of $p$ and for each $d$ prime to
  $p$ let $F_d=\Fq(t)$ with $t^d=u$.  Let $\rho$ be a representation
  satisfying the hypotheses of \ref{ss:rho-hyps} which is self dual of
  weight $w$ and sign $-1$.  Let $\n$ be the conductor of $\rho$,
  let $\n'$ be the part of $\n$ which is prime to the places $0$ and
  $\infty$ of $F$, and let $\swan_0(\rho)$ and $\swan_\infty(\rho)$ be
  the exponents of the Swan conductors of $\rho$ at $0$ and $\infty$.
  If $\deg\n'+\swan_0(\rho)+\swan_\infty(\rho)$ is odd then
  $\ord_{s=(w+1)/2}L(\rho,F_d,s)$ is unbounded as $d$ varies through
  integers prime to $p$.  More precisely, for $d$ of the form
  $d=q^n+1$ we have
\begin{equation*}
\ord_{s=(w+1)/2}L(\rho,F_d,s)\ge d/2n-c
\end{equation*}
and
\begin{equation*}
\ord_{s=(w+1)/2}L(\rho,\F_{q^{2n}}F_d,s)\ge d-c
\end{equation*}
where $c$ is a constant independent of $n$.
\end{thm}

\begin{proof}
  Clearly it suffices to prove the ``more precisely'' assertion.  The
  extension $F_d/F$ is geometrically abelian with geometric Galois
  group $H=\mu_d$ and the Frobenius $\phi$ acts on $\hat H=\Z/d\Z$ by
  multiplication by $q^{-1}$.  For $d$ of the form $q^n+1$ we have
  $q^n\equiv-1\pmod d$ and so every orbit $o$ of $C$ (the group
  generated by $\phi$) on $\Z/d\Z$ satisfies $o=-o$, i.e., is
  self-dual.  As pointed out in \ref{ss:special-d}, there are at least
  $(q^n-1)/2n$ self-dual orbits consisting of characters of higher
  order.

  Now for $n$ sufficiently large and all $o$ such that each $\chi\in
  o$ has sufficiently large order, the space of invariants of
  $\rho\tensor\chi$ under the inertia group at 0 or $\infty$ is
  trivial.  For such $n$ and $\chi$,
$$\deg\cond(\rho\tensor\chi)=
\deg(\n')+\swan_0(\rho)+\swan_\infty(\rho)+2\dim\rho$$ which is odd
and so Theorem~\ref{thm:Lzeroes} implies that for each such orbit $o$,
the $L$-function 
$L(\rho\tensor\sigma_{o,f},F,s)$ vanishes at $s=(w+1)/2$.  The number
of ``bad'' orbits is bounded independently of $n$ and so the
factorization in Subsection~\ref{ss:geom-abelian-base-change} shows
that $L(\rho,F_d,s)$ has a zero of order at least $d/2n-c$ at
$s=(w+1)/2$ for some constant $c$ independent of $n$.

Extending scalars to $\F_{q^{2n}}F$, each factor
$1-\left(Tq^{\frac{w+1}2}\right)^{|o|}$ dividing the $L$-function
becomes $\left(1-Tq^{2n\frac{w+1}2}\right)^{|a|}$ and so the total
order of vanishing of $L(\rho,\F_{q^{2n}}F(u^{1/d}),s)$ at $s=(w+1)/2$
is $\ge d-c$ for some constant $c$ independent of $n$.
\end{proof}

The analytic rank assertions in Theorems~\ref{thm:av1}, \ref{thm:av2},
and \ref{thm:ec1} will all be established using the ``towers''
Theorem~\ref{thm:towers}.  The Tate-Shafarevitch Theorem~\ref{thm:ec2}
will follow similarly from an orthogonal ($\sgn(\rho)=1$) variant of
the towers theorem.

We end this section with another example of towers leading to
unbounded ranks.  The proof is quite similar to that of
Theorem~\ref{thm:towers} and thus will be omitted.

\begin{thm}
  Let $E$ be an elliptic curve over a finite field $\Fq$ of
  characteristic $p$ and let $F=\Fq(E)$.  Let $\ell$ and $\ell'$ be
 \textup{(}not necessarily distinct\textup{)} prime numbers $\neq p$ with $\ell'$ odd.
  For each $n\ge1$ let $F_n=\Fq(E)$ and view $F_n$ as an extension of
  $F$ via pullback under the multiplication-by-$\ell^{\prime n}$
  isogeny $\ell^{\prime n}:E\to E$.  Assume that some power of
  Frobenius acting on the $\ell'$-torsion $E[\ell']$ has eigenvalue
  $-1$.  Let $\rho$ be an $\ell$-adic representation satisfying the
  hypotheses of \ref{ss:rho-hyps} which is self dual of weight $w$ and
  sign $-1$.  Assume that the degree of the conductor of $\rho$ is
  odd. Then
\begin{equation*}
\ord_{s=(w+1)/2}L(\rho,F_n,s)\ge n.
\end{equation*}
\end{thm}

\begin{rem}
  The hypothesis on $E$ in the theorem is very mild.  By assuming
  more about $E$ we can improve the lower bound in the theorem to
  $\Omega(\ell^{\prime n})$.  We omit the details.
\end{rem}

\section{Proof of the first part of Theorem~\ref{thm:av1}}\label{s:av1}
In this section we will show that there are many examples of curves
over $\Fp(t)$ whose Jacobians satisfy the first part of
Theorem~\ref{thm:av1}, namely they are absolutely irreducible,
non-isotrivial, and have large analytic rank.  To keep the exposition
brief, we have chosen examples where the necessary calculations have
already appeared in the literature, but the reader who is so inclined
will have no trouble finding many other examples. In the following two
sections we will give examples where one can also prove
the conjecture of Birch and Swinnerton-Dyer and therefore conclude
that algebraic ranks are also large.  As will be apparent, the class
of examples for which one can currently prove large algebraic ranks is
considerably smaller than that for which one can prove large analytic
ranks.

\subsection{}
Fix a prime $p$ and a positive integer $g$.  Let $F=\Fp(u)$ and for
each $d$ not divisible by $p$ let $F_d=\Fp(t)$ where $u=t^d$.  Suppose
that $C$ is a curve of genus $g$ smooth and proper over $F$, let
$J=J(C)$ be its Jacobian, and let $V=V_\ell J\tensor\Qlbar$ be the
$\ell$-adic Tate module of $J$ for some prime $\ell\neq p$.  Let
$\rho:G_F\to\aut(V^*)\cong\GL_{2g}(\Qlbar)$ be the
natural representation of Galois on $V^*\cong
H^1(C\times\overline{F},\Qlbar)$.  The representation $\rho$ satisfies
the hypotheses of Section~\ref{ss:rho-hyps} with weight $w=1$ and sign
$\sgn(\rho)=-1$.

The $L$-function of $J$ is of course the same as the $L$-function of
$\rho$ and so if $\rho$ satisfies the hypotheses of the towers
Theorem~\ref{thm:towers}, then $J$ will have large analytic rank over
$F_d$ for suitable $d$.

\subsection{}\label{ss:monodromy}
We consider the monodromy groups attached to $\rho$.  Let $\rho'$ be
the Tate twist $\rho\tensor \Qbar(1/2)$ which has weight $w=0$ and let
$\rho'_0$ be the restriction $\rho'|_{\gal(\overline{F}/\Fpbar F)}$.
Let $G^{\text{arith}}$ be the Zariski closure of the image of $\rho'$
and let $G^{\text{geom}}$ be the Zariski closure of the image of
$\rho'_0$.  The latter group is a (possibly non-connected) semi-simple
algebraic group over $\Qlbar$ and, because $\rho$ is self-dual of sign
$-1$, both groups are {\it a priori\/} contained in the symplectic
group $\Sp_{2g}$.  In the examples we will consider below, it will
turn out that $G^{\text{arith}}=G^{\text{geom}}=\Sp_{2g}$.

As usual, we say that $\rho'$ is irreducible if
$V^*$ has no non-trivial subspaces invariant under
$\rho'(G_F)$, or equivalently, under the action of
$G^{\text{arith}}$.  We say that $\rho'$ is Lie irreducible if the
restriction of $\rho'$ to any finite index subgroup of
$G_F$ is irreducible.  This is
equivalent to saying that $G^{\text{arith}}$ acts irreducibly and is
connected and in this case we also say that $G^{\text{arith}}$ acts
Lie irreducibly.

\subsection{}\label{ss:non-isotrivial}
We say that $J$ is non-isotrivial if there does not exist an abelian
variety $J_0$ defined over a finite field $\Fq$ and a finite extension
$K$ of $F=\Fp(t)$ containing $\Fq$ such that $J\times_FK\cong
J_0\times_{\Fq} K$.

It is clear that if the monodromy group $G^{\text{arith}}$ acts Lie
irreducibly, then $J$ is absolutely simple and non-isotrivial.
Indeed, if $J$ had a non-trivial isogeny decomposition over a finite
separable extension $K/F$, or if $J$ became isomorphic to a constant
abelian variety over a finite separable extension $K/F$, then $\rho'$
restricted to $G_K$ would be reducible.  Since the
monodromy group $G^{\text{arith}}$ is invariant under finite, purely
inseparable extensions similar statements hold for any finite
extension $K/F$.

In fact it is clear that when $G^{\text{arith}}$ acts Lie irreducibly
(as defined in \ref{ss:monodromy}),
$J$ is not even isogenous to a constant abelian variety over any
extension, since $G^{\text{arith}}$ is invariant under isogeny.
Therefore, for all finite extensions $K/F$, the $K/\Fq$-trace and
$K/\Fq$-image of $J\times_FK$ vanish.  (See Conrad~\cite{ConradK/k}
for a modern treatment of the $K/k$-trace and $K/k$-image.)

In light of this discussion, to prove the first part of
Theorem~\ref{thm:av1}, it will suffice to exhibit curves whose
Tate-module representations have $G^{\text{arith}}=\Sp_{2g}$ and which
satisfy the conductor hypothesis of the towers
Theorem~\ref{thm:towers}.  We do this in the following two
subsections.

\subsection{}\label{ss:analytic,p>2}
Assume that $p>2$ and choose a polynomial $f(x)\in\Fp[x]$ of degree
$2g$ with distinct roots, one of which is 0.  Consider the curve $C$
of genus $g$ smooth and proper over $F=\Fp(u)$ with affine equation
\begin{equation}\label{eqn:hyp-curve}
y^2=f(x)(x-u).
\end{equation}
Let $J$ be the Jacobian of $C$ and let $\rho$ be the representation of
$G_F$ on 
$H^1(C\times\overline{F},\Qlbar)$.  As we will explain in the rest of
this subsection, results of Katz and Sarnak show that $\rho$ satisfies
the hypotheses of Theorem~\ref{thm:towers} and has
$G^{\text{arith}}=\Sp_{2g}$ and so $J$ is an example satisfying the
first part of Theorem~\ref{thm:av1}.

In order to apply the results of Katz and Sarnak, we need to make one
translation.  Namely, they work with a lisse sheaf $\FF$ on an open
subset of $\P^1$ whose generic stalk is the cohomology with compact
supports of the affine curve defined by equation~\ref{eqn:hyp-curve}.
The smooth, proper model of this curve is obtained by adding exactly
one point at infinity, and so the compactly supported $H^1$ of the
open curve is canonically isomorphic to the usual $H^1$ of the proper
curve.  This implies that Katz and Sarnak's sheaf $\FF$ is the
restriction to an open of $\P^1$ of the middle extension sheaf
$\FF_\rho$ we considered in the proof of Theorem~\ref{thm:Lzeroes}.
The same issue arises in the next subsection, with the same
resolution.

Now by \cite{KS}*{10.1.12}, $\rho$ is everywhere tame and so
$\swan_0(\rho)=\swan_\infty(\rho)=0$.  By \cite{KS}*{10.1.9 and
  10.1.12}, $\n'$, the prime-to-zero-and-infinity part of the
conductor of $\rho$, is the sum of the zeros of $f$ except 0, each
taken with multiplicity one.  Thus, $\deg\n'=2g-1$ is odd and so
$\deg\n'+\swan_0(\rho)+\swan_\infty(\rho)$ is odd.  Thus
Theorem~\ref{thm:towers} implies that the analytic rank of $J$ is
unbounded in the tower of fields $F_d$.

By \cite{KS}*{10.1.16}, the geometric monodromy group of $\rho$ is the
full symplectic group $\Sp_{2g}$ and therefore the same is true of the
arithmetic monodromy group since $G^{\text{geom}}\subset
G^{\text{arith}}\subset\Sp_{2g}$.  Thus $J$ is absolutely simple and
non-isotrivial.

This completes the proof of the first part of Theorem~\ref{thm:av1}
for $p>2$.

\subsection{}\label{ss:analytic,p=2}
Now assume that $p=2$ and consider the curve $C$ of genus $g$ smooth
and proper over $\Fp(u)$ with affine equation
\begin{equation*}
y^2+xy=x^{2g+1}+ux.
\end{equation*}
Again let $\rho$ be the representation of $G_F$ on
$H^1(C\times\overline{F},\Qlbar)$.  It
is easy to see that $C$ has good reduction away from $u=0$ and
$\infty$ and so $\rho$ is unramified away from those two places.  By
\cite{KS}*{proof of 10.2.2}, $\rho$ is tamely ramified at $0$ and by
the first full paragraph of p.~302 of \cite{KS} and
\cite{KatzESDE}*{7.5.4}, the Swan conductor of $\rho$ at $\infty$ is
$2g-1$.  Thus $\deg\n'+\swan_0(\rho)+\swan_\infty(\rho)=2g-1$ is odd
and Theorem~\ref{thm:towers} shows that $J$ has unbounded
analytic rank in the tower of fields $F_d$.

By \cite{KS}*{10.2.2}, the geometric monodromy group of $\rho$ is
$\Sp_{2g}$ and so we conclude as in the previous section that $J$ is
absolutely simple and non-isotrivial.

This completes the proof of the first part of Theorem~\ref{thm:av1}
for $p=2$.

\section{BSD for curves defined by four monomials}
In this section we will show that the Jacobians of curves defined by
particularly simple equations satisfy the conjecture of Birch and
Swinnerton-Dyer.  The main tools are a beautiful observation of Shioda
\cite{ShiodaPicard}, already exploited in \cite{UlmerR}, that surfaces
defined by four monomials are often dominated by Fermat surfaces and
so satisfy the Tate conjecture, and well-known connections between the
conjectures of Tate and of Birch and Swinnerton-Dyer.

\subsection{}\label{ss:Shioda}
Let $k$ be a field and consider an irreducible polynomial $g\in
k[x_1,x_2,x_3]$ which is the sum of exactly 4 non-zero monomials:
$$g=c_0x_1^{a_{01}}x_2^{a_{02}}x_3^{a_{03}}+\cdots
+c_3x_1^{a_{31}}x_2^{a_{32}}x_3^{a_{33}}
=\sum_{i=0}^3c_i\prod_{j=1}^3x_j^{a_{ij}}$$ For $i=0,\dots,3$, let
$a_{i0}=1-\sum_{j=1}^3a_{ij}$ and let $A$ be the $4\times4$ integer
matrix $(a_{ij})$.  We say that {\it $g$ satisfies Shioda's
  conditions\/} if two requirements hold.  First, we require that the
determinant of $A$ be non-zero.  Assuming so, $A$ has an inverse in
$\GL_4(\Q)$ and there is a well-defined smallest positive integer
$\delta$ such that $B=\delta A^{-1}$ has integer coefficients.  Our
second requirement is that $\delta$ be non-zero in $k$.  Note that
$AB=\delta I_4$.  Note also that Shioda's conditions are independent
of the ordering of the variables $x_i$ and indeed independent of their
names, i.e., the condition makes sense for any polynomial ring in
three variables.  The reader may be surprised to see the 1 in the
definition of $A$ rather than $\deg(g)$, but it gives better (less
divisible) values of $\delta$.

\begin{thm}\label{thm:4-monomials}
  Let $F=\Fq(u)$ and let $X$ be a curve smooth and proper over $F$.
  Let $J$ be the Jacobian of $X$.  Assume that there exists an
  irreducible polynomial $g\in\Fq[u,x,y]\subset F[x,y]$ which is the
  sum of exactly 4 non-zero monomials, which satisfies Shioda's
  conditions, and which gives rise to the function field $F(X)$ in the
  following sense: $F(X)\cong\mathrm{Frac}(F[x,y]/(g))$.  Then the
  conjecture of Birch and Swinnerton-Dyer conjecture holds for $J$,
  namely $\rk J(F)=\ord_{s=1}L(J/F,s)$.
\end{thm}

\begin{rems}\hfil\break
\begin{enumerate}
\item \vskip-12pt It is known that over function fields the weak form
  of BSD ($\rk=\ord$) is equivalent to the finiteness of $\sha$ and
  implies the refined conjecture on the leading coefficient of the
  $L$-series.  We give a few more details about this below.
\item If $X$ satisfies the hypotheses of the theorem, then it is easy
  to check that the same is true for $X\times_F\F_{r}(t)$ where
  $t^d=u$ for any positive integer $d$ not divisible by $p$ and for
  any power $r$ of $q$.  Thus the theorem gives the truth of BSD for
  curves over towers of function fields.
\end{enumerate}
\end{rems}

\begin{proof}
[Proof of Theorem~\ref{thm:4-monomials}]
Let $\ZZ$ be the surface in $\A^3$ defined by $g=0$.  Since $g$ is
irreducible, $\ZZ$ is reduced and irreducible and so has a dense open
subset smooth over $\Fq$.  There is a morphism $\ZZ\to\A^1$, namely
$(u,x,y)\mapsto u$.

Let $\XX$ be a model of $\ZZ$ smooth and proper over $\Fq$.  There is
a rational map $\XX\ratto\P^1$ and at the expense of blowing up and
down we may assume that we have a morphism $\pi:\XX\to\P^1$ which is
relatively minimal.  The generic fiber of $\pi$ is a regular scheme
(since its local rings are local rings of $\XX$) of dimension 1 which
is proper over $F$ and has the same function field as $X$ does and so
it is isomorphic to $X\to\spec F$.  


Let $F_\delta$ be the Fermat surface of degree $\delta$.  The
assumption that $g$ satisfies Shioda's conditions implies that
(extending the ground field $\Fq$ if necessary) there is a dominant
rational map $\FF_\delta\ratto\ZZ$ and therefore also a dominant rational
map $\FF_\delta\ratto\XX$.  Indeed, let us use the notation of
Subsection~\ref{ss:Shioda} (setting $u=x_1$, $x=x_2$, and $y=x_3$) and
define $\ZZ'$ as the zero set in $\P^3\setminus\{x_0=0\}$ of the
homogeneous Laurent polynomial $\sum_ic_i\prod_jx_j^{a_{ij}}$.
Clearly $\ZZ'$ is birational to $\ZZ$.  If $y_0,\dots,y_3$ are the
standard coordinates on $F_\delta$ (so that
$y_0^\delta+\cdots+y_3^\delta=0$), then a rational map $F_\delta\ratto
\ZZ'$ is given by $x_j\mapsto\prod_kd_ky_k^{b_{jk}}$ where
$(b_{jk})=B=\delta A^{-1}$ and $d_0,\dots,d_3\in\Fqbar$ are solutions
to $\prod_jd_j^{a_{ij}}=c_i^{-1}$.  This proves that there is a
dominant rational map from $F_\delta$ to $\XX$.

Now the Fermat surface $\FF_\delta$ is dominated by a product of
Fermat curves \cite{ShiodaKatsura} and therefore so is $\XX$.  As we
will explain presently, this domination by a product of curves is
enough to imply the Tate conjecture for $\XX$ and the conjecture of
Birch and Swinnerton-Dyer for $J$.  (There is a large literature on
the connection between these conjectures; the approach that follows is
perhaps ahistorical, but has the virtue of being efficient and
clear-cut.)  That $\XX$ is dominated over an extension of $\Fq$ by a
product of curves implies \cite{TateM}*{\S5} that the Tate conjecture
(relating the rank of the N\'eron-Severi group of $\XX$ to it zeta
function) holds for $\XX$ over an extension of $\Fq$ and therefore
also over $\Fq$.  Consideration of the Kummer sequence in \'etale
cohomology \cite{TateB}*{5.2} implies that the $\ell$-primary part of
the Brauer group of $\XX$ is finite for all $\ell\neq p$.  A theorem
of Artin generalized by Grothendieck implies that the same holds for
the $\ell$-primary parts of the Tate-Shafarevitch group of $J$ over
$F$.  (This can be extracted from \cite{GrothB3}*{\S4}; in the case
when $X$ has an $F$-rational point, it is proven in
\cite{GrothB3}*{\S4} that $Br(\XX)=\sha(J/F)$.)  Finally, a recent
paper of Kato and Trihan \cite{KatoTrihan} proves that the full
conjecture of Birch and Swinnerton-Dyer holds for $J$ over $F$ as soon
as one $\ell$-primary part of $\sha$ is finite.  (In the applications
below, we will only need the theorem in cases where $X$ is a curve
with an $F$-rational point and in this case, the reference to
\cite{KatoTrihan} may be replaced with the simpler
\cite{MilneADT}*{III.9.7}.)
\end{proof}

\begin{rem}
  The above proof of the Tate conjecture for $\XX$ ultimately comes
  down to two facts: Tate's theorem on isogenies of abelian varieties
  over finite fields and the fact that $\XX$ is dominated by a product
  of curves, Fermat curves as it turns out.  It is not difficult to
  produce examples of surfaces not defined by four monomials which are
  dominated by products of curves.  But the four monomials property
  has the charm that it is obviously preserved in towers, i.e., when
  $u$ is replaced by $t^d$.  It looks like an interesting problem to
  give examples of towers of surfaces dominated by products of curves
  beyond the four monomial case.
\end{rem}

\section{End of the proof of Theorem~\ref{thm:av1}}\label{s:algranks}
In order to finish the proof of Theorem~\ref{thm:av1} we have to
exhibit for every prime $p$ and every integer $g>0$ a curve $X$ smooth
and proper over $F=\Fp(u)$ of genus $g$ with three properties: (i) the
Galois representation $\rho$ on $H^1(X\times\overline{F},\Qlbar)$ should
satisfy the hypotheses of the towers Theorem~\ref{thm:towers} (so that
we have large analytic ranks); (ii) the monodromy group of $\rho$
should be Lie-irreducible (so that $J=J(X)$ will be absolutely simple
and non-isotrivial, cf.~Subsection~\ref{ss:non-isotrivial}); and (iii)
$X$ should satisfy the hypotheses of the ``four monomials''
Theorem~\ref{thm:4-monomials} (so that the conjecture of Birch and
Swinnerton-Dyer holds for $J$).  In this section we will exhibit
curves with these properties.

\subsection{}
There are many examples of such curves.  The ones we have chosen here
allow for a fairly unified treatment and, as will be explained in
Section~\ref{s:RankBounds}, they have good properties with respect to
rank bounds.  Here are the curves we will study:
\begin{align}
\label{eqn:case1}&p>2\quad p\nodiv(2g+2)(2g+1)&&y^2=x^{2g+2}+x^{2g+1}+u\\
\label{eqn:case2}&p>2\quad p|(2g+2)&&y^2=x^{2g+2}+x^{2g+1}+ux\\
\label{eqn:case3}&p>2\quad p|(2g+1)&&y^2=x^{2g+1}+x^{2g}+ux\\
\label{eqn:case4}&p=2&&y^2+xy=x^{2g+1}+ux
\end{align}

\subsection{}\label{ss:4monsok}
Let $F=\Fp(u)$ and let $X$ be the regular, proper model of one of the
affine curves over $F$ defined by equations
\ref{eqn:case1}-\ref{eqn:case4}.  Then it is easy to see that $X$ is
smooth over $\spec F$ and satisfies the hypotheses of the four
monomials Theorem~\ref{thm:4-monomials}.  (The quantity $\delta$
appearing in Shioda's conditions is equal to 2 in the first three
cases, and $2g+1$ in the fourth case.)  Thus the Jacobian $J$ of the
curve $X$ and all its base changes under $u\mapsto t^d$ satisfy the
conjecture of Birch and Swinnerton-Dyer.

In Subsections~\ref{lemma:tameness} to \ref{lemma:LieIrreducible}
below, we will prove that the curve defined by
equation~\ref{eqn:case1} satisfies the hypotheses of the towers
Theorem~\ref{thm:towers} and has Lie irreducible monodromy.  Then in
Section~\ref{ss:modifications} we will give the minor modifications
needed to treat equations~\ref{eqn:case2} and \ref{eqn:case3}.
Finally, in Section~\ref{ss:p=2alg} we treat the last case,
equation~\ref{eqn:case4}.

\begin{lemma}\label{lemma:tameness}
  Suppose that $F$ is a global field of characteristic $p>2$ and $X$
  is a smooth hyperelliptic curve over $F$ with affine equation
  $y^2=f(x)$ for some polynomial $f\in F[x]$ with distinct roots.  Let
  $\rho$ be the natural representation of $G_F$ on
  $H^1(X\times\overline{F},\Qlbar)$ ($\ell$ any prime $\neq p$).  Then
  $\rho$ is everywhere tamely ramified if and only if the splitting
  field of $f$ is an everywhere tamely ramified extension of $F$.
\end{lemma}

\begin{proof}
  The question of whether or not $\rho$ is everywhere tame is
  independent of $\ell$ by \cite{Saito}*{3.11} and thus we may assume
  $\ell=2$.  As a representation of $G_F$,
  $H^1(X\times\overline{F},\Q_2)$ is dual to the 2-adic Tate module
  $V_2J$ and so it suffices to show that the latter is everywhere
  tame.  The 2-torsion in $J$ is spanned by the classes of divisors of
  degree zero supported on the Weierstrass points $(\alpha,0)$ where
  $\alpha$ is a root of $f$ together with the point at infinity on $X$
  if the degree of $f$ is odd.  Using this it is not hard to check
  that the fixed field of the kernel of the action of Galois on the
  2-torsion $J[2]$ is precisely the splitting field of $f$.  The
  restriction of $\rho$ to the Galois group of this field takes its
  values in $I+2M_2(\Z_2)$ which is a pro-2 group and so this
  restriction is at worst tamely ramified.  Therefore $\rho$ itself is
  tamely ramified if and only if the splitting field of $f$ is tamely
  ramified over $F$.
\end{proof}

\begin{lemma}\label{lemma:Sn}
  Let $p>2$ be a prime and let $n>1$ be an integer such that $p\nodiv
  n(n-1)$.  Then $f(x)=x^n+x^{n-1}+u$ is irreducible over $\Fpbar(u)$
  and its splitting field is everywhere tame with Galois group $S_n$,
  the symmetric group on $n$ letters.
\end{lemma}

\begin{proof}
  We will see below that $f$ is irreducible over $\Fpbar(u)$.  Let $K$
  be its splitting field.  Considering $f$ and its derivative
  $f'=nx^{n-1}+(n-1)x^{n-2}$ we see that the reduction of $f$ at a
  place of $\Fpbar(u)$ has distinct roots except at the places $u=0$,
  $u=a:=-(1-n)^{n-1}n^{-n}$, and $u=\infty$.  Thus $K$ is unramified
  over $\Fpbar(u)$ away from these places.

  Let $F_0$ be the completion $\Fpbar((u))$ of $\Fpbar(u)$ at $u=0$.
  Consideration of the Newton polygon of $f$ with respect to the
  valuation of $F_0$ shows that one root of $f$, call it $\alpha$,
  lies in $F_0$ and is congruent to $-1$ modulo $u$ and the other
  roots have valuation $1/(n-1)$.  We have
$$f(x)=(x-\alpha)\left(x^{n-1}+(\alpha+1)x^{n-2}+\cdots
  +(\alpha^{n-2}+\alpha^{n-3})x-u/\alpha\right).$$ 
If $g(x)$ denotes the second factor on the right and $\beta$ is a root
of $g$, then $g(x/\beta)$ is congruent modulo the maximal ideal of
$F_0(\beta)$ to $x^{n-1}+b$ where $b\neq0$.  Since $p\nodiv(n-1)$,
this reduction has distinct roots and so by Hensel's lemma all the
roots of $g$ lie in $F_0(\beta)$.  Thus the splitting field of $f$
over $F_0$ is a cyclic extension of degree $(n-1)$.  We conclude that
there are two places of $K$ over $u=0$, one unramified and the other
totally ramified with index $e=n-1$.  In particular, the ramification
over 0 is tame.

Now let $v=1/u$ and $F_\infty=\Fp((v))$ be the completion of
$\Fpbar(u)$ at $u=\infty$.  Changing variables, let
$$h(x)=u^{-n}f(ux)=x^n+vx^{n-1}+v^{n-1}.$$
Consideration of the Newton polygon of $h$ with respect to the
valuation of $F_\infty$ shows that $h$ is irreducible over $F_\infty$
(and so $f$ is irreducible over $\Fpbar(u)$) and that its roots all
have valuation $(n-1)/n$.  If $\gamma$ is one of these roots, then
modulo the maximal ideal of $F_\infty(\gamma)$, $h(x/\gamma)$ is
congruent to $x^n+c$ with $c\neq0$.  Since $p\nodiv n$ this reduction
has distinct roots and so by Hensel's lemma, all of the roots of $h$
lie in $F_\infty(\beta)$.  Therefore the splitting field of $f$ over
$\F_\infty$ is a cyclic extension of degree $n$ and is totally
ramified ($e=n$).  Thus $K$ is totally and tamely ramified over the
place $u=\infty$ of $\Fpbar(u)$.

Now let $v=u-a$ and $F_a=\Fpbar((v))$ be the completion of $\Fpbar(u)$
at $u=a$.  The specialization of $f$ to $u=a$ has $n-2$ simple roots
and a double root $b=(1-n)/n$.  By Hensel's lemma, $f$ has $n-2$ of
its roots in $F_a$.  Considering the Newton polygon of $f(x+b)$, we
see that the other two roots of $f$ have valuation $1/2$ and thus lie
in a ramified quadratic extension of $F_a$.  Thus $K$ has $n-2$ split
places ($e=1$) and one place with $e=2$ over $u=a$.  Since $p>2$, the
ramification is tame.

This shows that $K$ is everywhere tame over $\Fpbar(u)$.  Inertia at 0
is generated by an $(n-1)$-cycle $\sigma$, inertia at $\infty$ is
generated by an $n$-cycle $\tau$ and inertia at $a$ is generated by a
simple transposition $\rho$.  Moreover, by the known structure of the
tame fundamental group of $\P^1$ minus 3 points, we may choose these
generators so that $\rho\sigma=\tau$.  Choosing labels so that 1 is
the fixed point of $\sigma$ and $\tau=(12\cdots n)$ we see immediately
that $\rho=(12)$.  Since the symmetric group is generated by $(12)$
and $(12\cdots n)$ we conclude that the Galois group of $K$ over
$\Fpbar(u)$ is $S_n$.
\end{proof}

\begin{cor}\label{cor:2components}
  With hypotheses as in Lemma~\ref{lemma:Sn}, The affine plane curve
  defined by
$$g(x,x')=x^n+x^{n-1}-x^{\prime n}-x^{\prime n-1}=0$$ 
has exactly two irreducible components over $\Fpbar$ both of which are
rational over $\Fp$.
\end{cor}

\begin{proof}
  With notation as in Lemma~\ref{lemma:Sn}, the curve in question is
  the fiber product of two copies of $f=0$ over the $u$-line.  Its set
  of irreducible components is thus in bijection with the orbits of
  $S_n$ on the set of ordered pairs of roots of $f$ in
  $\overline{\Fp(u)}$ and there are two such orbits, the diagonal and
  the rest.  The equations of the two components are $x-x'=0$ and
  $g(x,x')/(x-x')=0$, both of which are $\Fp$-rational.
\end{proof}

\subsection{}
From here through \ref{lemma:LieIrreducible} we let $X$ be the curve
defined by \ref{eqn:case1} and $J$ its Jacobian.  We let $\rho$ be the
representation of $G_F$ on
$H^1(X\times\overline{F},\Qlbar)$.  Let $\pi:\XX\to\P^1$ be the model of
$X$ constructed in the proof of Theorem~\ref{thm:4-monomials}.  Then
$j:U\into\P^1$, with $U=\P^1\setminus\{0,a,\infty\}$, is the largest
open subset of $\P^1$ over which $\pi$ is smooth.  Let $\FF_U$ be the
lisse sheaf on $U$ corresponding to the representation $\rho$ and set
$\FF=j_*\FF_U$.  This is the ``middle extension'' sheaf attached to
$\rho$ and we may recover $\rho$ from it as the stalk
$\FF_{\overline{\eta}}$ at the geometric generic point $\overline\eta$
corresponding to the fixed algebraic closure $F^{\rm alg}$ of $F$.


Recall that a linear transformation is called a unipotent
pseudoreflection if all of its eigenvalues are 1 and its space of
invariants has codimension 1.

\begin{lemma}\label{lemma:Lefschetz}
  The sheaf $\FF$ is everywhere tamely ramified.  At the place $u=a$
  inertia acts via unipotent pseudoreflections and in particular the
  exponent of the Artin conductor is 1.  Therefore, $\FF$ \textup{(}or
  rather $\rho$\textup{)} satisfies the hypotheses of
  Theorem~\ref{thm:towers}.
\end{lemma}

\begin{proof}
  The preceding two lemmas show that $\rho$ is everywhere tamely
  ramified and therefore the same is true of $\FF$.  To analyze the
  ramification at $u=a$, consider the following surface: let
\begin{align*}
V_1&=\spec\Fp[u,x,y]/\left(y^2-(x^{2g+2}+x^{2g+1}+u)\right)\cr
V_2&=\spec\Fp[u,x',y']/\left(y^{\prime2}-(1+x'+ux^{'2g+2})\right)
\end{align*}
and define $\YY$ as the result of glueing $V_1$ and $V_2$ using the
map
$$(x',y',u)=(x^{-1},yx^{-g-1},u).$$  
There is a map $\YY\to\A^1$ (projection onto the $u$ coordinate) which
is proper, relatively minimal, and whose generic fiber is $X$.
Moreover, $\YY$ is a regular surface and therefore we may identify $Y$
with the open subset $\pi^{-1}(\A^1)\subset\XX$ where $\pi:\XX\to\P^1$
was constructed in the proof of Theorem~\ref{thm:4-monomials}.  The
restriction $\pi_{\A^1}:\YY\to\A^1$ is smooth except over $u=a$, and
over $u=a$, it has an isolated singularity which is an ordinary double
point.  In classical language, $\pi_{\A^1}$ is a Lefschetz pencil.
The famous Picard-Lefschetz formula (\cite{SGA7-2}*{3.4} or
\cite{MilneEC}*{V.3.14}) gives the action inertia on
$\FF_{\A^1}=R^1\pi_{\A^1*}\Qlbar$ and in particular it shows that the
action of inertia at $u=a$ is by unipotent pseudoreflections.  Since
the ramification at $u=a$ is tame, the exponent of the Artin conductor
is just the codimension of the space of inertia invariants which is 1.
\end{proof}

\begin{lemma}\label{lemma:irreducible}
  $\FF$ is geometrically irreducible as middle extension sheaf on
  $\overline{\P^1}$.  Equivalently, $\rho$ restricted to
  $\gal(\overline{F}/\Fpbar F)$ is geometrically irreducible.
\end{lemma}

\begin{proof}
  We apply the diophantine criterion for irreducibility of Katz, along
  the lines of \cite{KS}*{10.1.15}.  This criterion amounts to using
  the Grothendieck-Lefschetz trace formula and the Weil conjectures to
  prove that $\en_{\gal(\overline{F}/\Fpbar F)}(\rho)$ is
  one-dimensional by estimating certain sums of traces.  By Schur's
  lemma, this one-dimensionality is equivalent to the desired
  irreducibility.  We refer to \cite{KS}*{10.1.15} for more details.

What has to be shown is that 
$$\sum_{\substack{u\in\F_q\\u\neq0,a}}Tr(Fr_u|\FF)^2=q^2+O(q^{3/2})$$
as $q$ ranges through all powers of $p$.

Using the Grothendieck-Lefschetz trace formula on the fibers of $\pi$
we see that
\begin{align*}
Tr(Fr_u|\FF)&=q-1-\sum_{x\in\Fq}\left(1+\chi(x^{2g+2}+x^{2g+1}+u)\right)\\
&=-1-\sum_{x\in\Fq}\chi(x^{2g+2}+x^{2+1}+u)
\end{align*}
where $\chi$ is the nontrivial quadratic character of $\Fq^\times$
extended as usual to a function on $\Fq$. (The reader will note that
there are two points at infinity on the affine curve \ref{eqn:case1}.)
By theorems of Weil and Deligne, the trace is $O(q^{1/2})$.  Thus the
sum to be estimated is
\begin{equation*}
\sum_{\substack{u\in\F_q\\u\neq0,a}}
\left(-1-\sum_{x\in\Fq}\chi(x^{2g+2}+x^{2g+1}+u)\right)^2
\end{equation*}
and because of the Deligne estimate, we may drop the conditions
$u\neq0,a$.  Thus our sum is
\begin{multline*}
\sum_{u\in\F_q}\sum_{x,x'\in\Fq}
\chi\left((x^{2g+2}+x^{2g+1}+u)(x^{\prime2g+2}+x^{\prime2g+1}+u)\right)+O(q^{3/2})\\
=\sum_{x,x'\in\Fq}\sum_{u\in\Fq}
\chi\left(u^2+u(x^{2g+2}+x^{2g+1}+x^{\prime2g+2}+x^{\prime2g+1})\right.\\
+\left.(x^{2g+2}+x^{2g+1})(x^{\prime2g+2}+x^{\prime2g+1})\right)+O(q^{3/2}).
\end{multline*}
The inner sum over $u$ is related to the number of points on the
hyperelliptic curve
$$y^2=u^2+u(x^{2g+2}+x^{2g+1}+x^{\prime2g+2}+x^{\prime2g+1})+
(x^{2g+2}+x^{2g+1})(x^{\prime2g+2}+x^{\prime2g+1}).$$ 
Noting that
$g(x,x')=(x^{2g+2}+x^{2g+1})-(x^{\prime2g+2}+x^{\prime2g+1})=0$ if and
only if the quadratic polynomial in $u$ has a double root, we see that
the sum over $u$ is $-1$ if $g(x,x')\neq0$ and it is $q-1$ if
$g(x,x')=0$.  Therefore the sum to be estimated is
$$\sum_{\substack{x,x'\in\Fq\\g(x,x')\neq0}}(-1)+
\sum_{\substack{x,x'\in\Fq\\g(x,x')=0}}(q-1)+O(q^{3/2}).$$ 
The first sum is over an affine open subset of $\A^2$ and is therefore
$-q^2+O(q)$.  The second sum is over a curve which by
Corollary~\ref{cor:2components} has exactly 2 components and therefore
has $2q+O(q^{1/2})$ points.  Thus the second sum contributes
$2q^2+O(q^{3/2})$ and the entire sum is $q^2+O(q^{3/2})$.
\end{proof}

\begin{rem}
  Another approach to irreducibility would be to assume $\ell=2$ and
  use Lemma~\ref{lemma:Sn} to argue that the mod 2 representation
  $J[2]$ is irreducible and so {\it a fortiori\/} the $2$-adic
  representation $V_2J$ is irreducible.  But this argument does not
  apply to the other curves \ref{eqn:case2} and \ref{eqn:case3}
  whereas the argument given above does apply with minor
  modifications.
\end{rem}

We recall from Subsection~\ref{ss:monodromy} the notion of a Lie
irreducible representation.

\begin{lemma}\label{lemma:LieIrreducible}
  The representation $\rho$ is Lie irreducible.  In fact, the
  geometric monodromy group of $\rho$ is the full symplectic group
  $\Sp_{2g}$.
\end{lemma}

\begin{proof}
  Later in the proof we are going to assume that $\ell=2$.  It follows
  from \cite{Chin} and the fact that $G^\text{geom}$ is {\it a
    priori\/} contained in $\Sp_{2g}$ that if $G^\text{geom}=\Sp_{2g}$
  for one $\ell\neq p$ then it is so for all $\ell\neq p$.  The reader
  who prefers not to go into this may simply assume that $\ell=2$ for
  the entire proof of the last part of Theorem~\ref{thm:av1}.

  We have already seen that $\rho':=\rho|_{\gal(\overline{F}/\Fpbar
    F)}$ is irreducible.  It follows from \cite{KatzMG}*{Proposition
    1} that $\rho'$ is either Lie irreducible, tensor decomposable (in
  the sense that $\rho'\cong\sigma\tensor\tau$ where $\sigma$ and
  $\tau$ are representations of $\gal(\overline{F}/\Fpbar F)$ with
  $\sigma$ Lie irreducible and $\tau$ of degree $>1$ with finite
  image), or is induced from a representation of a proper subgroup of
  $\gal(\overline{F}/\Fpbar F)$.  We rule out the last two
  possibilities.

  Since inertia at $u=1$ acts via a unipotent pseudoreflection,
  $\rho'$ can be tensor decomposable only if the degree of $\sigma$ is
  $1$, which implies that $\sigma$ has finite image.  But if it were
  so, then $\rho'$ would also have finite image and this is impossible
  because a unipotent pseudoreflection has infinite order.

  Now suppose that $\rho'$ is induced from a representation $\sigma$
  of some subgroup of $\gal(\overline{F}/\Fpbar F)$ corresponding to a
  cover $\pi:\Curve\to\overline{\P^1}$.  We write $\rho'(x)$ for
  $\rho$ viewed (by restriction) as a representation of the inertia
  group $I(x)$ at a place $x\in\overline{\P^1}$ and similarly for
  $\sigma$.  We have
$$\rho'(x)=\bigoplus_{y\mapsto x}\ind^{I(x)}_{I(y)}\sigma(y)$$
where the sum is over points of $\Curve$ mapping to $x$.  Considering
this equality at the place $u=a$ where inertia acts with invariants of
codimension 1, we see that $I(y)=I(a)$ for all $y$ over $a$, in other
words, $\sigma$ must be unramified over $a$.  We also must have that
$\sigma$ is at worst tamely ramified over $0$ and $\infty$.  This
means that $\pi$ must be a cyclic cover obtained by extracting a root
of $u$, i.e., $\Curve=\P^1$ with coordinate $v$ and $\pi^*(u)=v^m$.

We argue that $m>1$ is incompatible with what we know about the action
of inertia at $u=0$ on 2-torsion.  Indeed, if $m>1$ and $h$ is a
generator of tame inertia at $0$, then the trace of $h$ on $\rho'(0)$
is zero.  On the other hand, the action of $h$ on the 2-torsion
subgroup $J[2]$ has non-zero trace.  More precisely, label the roots
$\alpha_i$ ($i=0,\dots,2g+1$) of $f(x)=x^{2g+2}+x^{2g+1}+t$ so that
$\alpha_0$ is fixed by $I(0)$ and the others are permuted cyclically
(cf.~the proof of Lemma~\ref{lemma:Sn}).  Let $P_i$ be the Weierstrass
point $(\alpha_i,0)$ on $X$ and let $D_i$ be the class of $P_i-P_0$ in
$J[2]$.  Then $J[2]$ is generated as $\F_2$ vector space by
$D_1,\dots,D_{2g+1}$ with the single relation $\sum D_i=0$.  The proof
of \ref{lemma:Sn} shows that $h$ permutes the $D_i$ cyclically.
Considering the matrix of $h$ in the basis $D_1,\dots,D_{2g}$ of
$J[2]$ we see that $h$ has trace $-1=1$ on $J[2]$ and therefore $h$
has non-zero trace on $\rho'(0)$.  Thus $m>1$ is impossible and so
$\FF$ is not induced from any non-trivial cover.

This completes the proof that $\rho$ is Lie irreducible.  Since there
is a place where inertia acts on $\rho$ by a unipotent
pseudoreflection, in fact the geometric monodromy group of $\rho$ is
$Sp_{2g}$ (see \cite{KatzESDE}*{1.5} for details).
\end{proof}

\subsection{}\label{ss:modifications}
This completes the proof that the curves \ref{eqn:case1} are examples
for Theorem~\ref{thm:av1}.  For \ref{eqn:case2} and \ref{eqn:case3},
the proof is essentially the same.  Very minor modifications are
needed in Lemma~\ref{lemma:Lefschetz} (checking that $\XX\to\P^1$ is a
Lefschetz pencil over $\A^1\setminus\{0\}$) and in
Lemma~\ref{lemma:irreducible} (applying the diophantine criterion for
irreducibility).  At the end of the proof of
Lemma~\ref{lemma:LieIrreducible} (checking Lie irreducibility), for
\ref{eqn:case2}, we use the action of monodromy at $\infty$ (which
again is a cyclic permutation of order $2g+1$) and for \ref{eqn:case3}
we use Lemma~\ref{lemma:Sn} to see that in terms of a suitable basis
of $J[2]$, a generator $h$ of tame inertia at 0 has one fixed vector
and permutes the other $2g-1$ vectors cyclically.  Altogether this
completes the proof of Theorem~\ref{thm:av1} for $p>2$.

\subsection{}\label{ss:p=2alg}
For $p=2$ and the curve \ref{eqn:case4}, we already saw in
Subsection~\ref{ss:analytic,p=2} that the hypotheses of the towers
Theorem~\ref{thm:towers} are satisfied and that this curve has large
geometric monodromy group.  We also saw in Subsection~\ref{ss:4monsok}
that the hypotheses of the four monomials
Theorem~\ref{thm:4-monomials} are satisfied.  This completes the proof
of Theorem~\ref{thm:av1} for $p=2$.

\section{Proof of Theorem~\ref{thm:av2}}\label{s:av2}
We will use some basic representation theory to see that the action of
$G_F$ on $H^k_{prim}(J,\Qlbar)$ satisfies the hypotheses
of Theorem~\ref{thm:towers} where $J$ is the Jacobian of one of the
curves studied in Section~\ref{ss:analytic,p>2} for suitable $g$ and
$k$.

\subsection{}\label{ss:la}
If $A$ is a principally polarized abelian variety of dimension $g$
over a field $F$ (e.g., a Jacobian), then the \'etale cohomology group
$H^1(A\times\overline{F},\Qlbar)$ carries a non-degenerate symplectic
form $\langle ,\rangle $.  If $k$ is odd, then the same is true of
$H^k(A\times\overline{F},\Qlbar)
\cong\bigwedge^kH^1(A\times\overline{F},\Qlbar)$.

If $k>1$, define a linear map
$\phi:\bigwedge^kH^1(A\times\overline{F},\Qlbar)
\to\bigwedge^{k-2}H^1(A\times\overline{F},\Qlbar)$
by
\begin{equation*}
\phi(w_1\wedge\cdots\wedge w_k)
=\sum_{1\le i<j\le k}(-1)^{i+j-1}
\langle w_i,w_j\rangle w_1\wedge\cdots\wedge \hat w_i\wedge
\cdots\wedge\hat w_j\wedge\cdots\wedge w_k
\end{equation*}
where a hat denotes a vector to omit.  

Choose a symplectic basis $v_1,\dots,v_{2g}$ of
$H^1(A\times\overline{F},\Qlbar)$, i.e., one satisfying $\langle
v_i,v_{j+g}\rangle =\delta_{ij}$ and $\langle v_i,v_j\rangle=\langle
v_{i+g},v_{j+g}\rangle=0$ for $1\le i,j\le g$.  Define a linear map
$\psi:\bigwedge^{k-2}H^1(A\times\overline{F},\Qlbar)
\to\bigwedge^{k}H^1(A\times\overline{F},\Qlbar)$
by
\begin{equation*}
\psi(w)=w\wedge(v_1\wedge v_{g+1} + \cdots+v_g\wedge v_{2g}).
\end{equation*}
The map $\psi$ is independent of the choice of basis; it depends only
on the symplectic form $\langle ,\rangle $.

We define the primitive part $H^k_{prim}(A\times\overline{F},\Qlbar)$ to
be the kernel of $\phi$.  It is well known that if $k\le g$ then
$\psi$ is injective and we have a direct sum decomposition
\begin{equation}\label{eqn:decomp}
H^k(A\times\overline{F},\Qlbar)\cong
H^k_{prim}(A\times\overline{F},\Qlbar)\bigoplus \im(\psi).
\end{equation}
Also, both $\phi$ and $\psi$ are equivariant for the action of
$G_F$ and so this direct sum decomposition is
respected by Galois.  The symplectic pairing on
$H^k(A\times\overline{F},\Qlbar)$ induced by $\langle ,\rangle $
restricts to a non-degenerate symplectic pairing on
$H^k_{prim}(A\times\overline{F},\Qlbar)$.

\subsection{}
Now assume that $p>2$ and let $C$ be one of the curves considered in
Section~\ref{ss:analytic,p>2} and $J$ its Jacobian.  Let $\rho_1$ be
the representation of $G_F$ on
$H^1(J\times\overline{F},\Qlbar)$ and let $\rho_k$ be the representation
of $G_F$ on $H^k_{prim}(J\times\overline{F},\Qlbar)$.
Then $\rho_k$ is self-dual of weight $k$ and sign $-1$.  Because the
monodromy group of $\rho_1$ is the full symplectic group, $\rho_k$ is
irreducible.

The representation $\rho_1$ is everywhere tamely ramified and at each
finite place $u_0\in\A^1$ of bad reduction, the local inertia group
$I(u_0)$ acts via unipotent pseudoreflections.  (By
\cite{KS}*{10.1.13} the action is either trivial or by unipotent
pseudoreflections, but it is not hard to see that the fibration
$\XX\to\P^1$ attached to the curves we are considering is a Lefschetz
fibration over $\A^1$ and so the local monodromies at finite places of
bad reduction are in fact unipotent pseudoreflections.)  This means
that in terms of a suitable symplectic basis $v_1,\dots,v_{2g}$, the
action of inertia is
\begin{align}\label{eqn:I-action}
s(v_1)&=v_1+\lambda(s)v_{g+1}\cr
s(v_i)&=v_i\qquad i=2,\dots,2g
\end{align}
where $\lambda:I(u_0)\to\Qlbar$ is a non-zero character.

Because $\rho_1$ is everywhere tame, $\rho_k$ is also everywhere tame.
In particular, the exponent of the conductor of $\rho_k$ at a place
$u_0$ is just the codimension of the space of invariants of
$\rho_k(I(u_0))$.  We will show that at finite $u_0$ this codimension
depends only on $g$ and $k$ and that for every $k$ there are
infinitely many $g$ such that this codimension is odd.  For such $g$
and $k$, the representation $\rho_k$ satisfies the hypotheses of
Theorem~\ref{thm:towers} and this will prove Theorem~\ref{thm:av2}.

\subsection{}
We fix a finite place $u_0\neq0$ of $F$ where $C$ has bad reduction
and consider the action of the inertia group $I(u_0)$ on
$\bigwedge^kH^1(J\times\overline{F},\Qlbar)$.  Choose a symplectic basis
$v_1,\dots,v_{2g}$ of $V=H^1(A\times\overline{F},\Qlbar)$ such that the
action of $I(u_0)$ is given by \ref{eqn:I-action}.  Let $V_1$ be the
span of $v_1$ and $v_{g+1}$ and let $V_2$ be the span of
$v_2,\dots,v_g$ and $v_{g+2},\dots,v_{2g}$.  Then we have
\begin{equation*}
\bigwedge^kV\cong(\bigwedge^2V_1\tensor\bigwedge^{k-2}V_2)\bigoplus
(V_1\tensor\bigwedge^{k-1}V_2)\bigoplus
(\bigwedge^{k}V_2).
\end{equation*}
It is easy to see using \ref{eqn:I-action} that $I(u_0)$ acts
trivially on the first and third summands and that the codimension of
its invariants on the middle summand is
$\dim\bigwedge^{k-1}V_2=\binom{2g-2}{k-1}$.

A similar analysis shows that the space of $I(u_0)$-invariants on
$\psi(\bigwedge^{k-2}V)$ has codimension $\binom{2g-2}{k-3}$.  From
the direct sum decomposition \ref{eqn:decomp} we conclude that the
exponent of the conductor of $\rho_k$ at $u_0$, i.e., the codimension
of the $I(u_0)$ invariants on $H^k_{prim}(J\times\overline{F},\Qlbar)$,
is $\binom{2g-2}{k-1}-\binom{2g-2}{k-3}$.

\subsection{}
To conclude the proof of Theorem~\ref{thm:av2}, we will show that for
every odd $k$ there are infinitely many $g$ such that
$\binom{2g-2}{k-1}-\binom{2g-2}{k-3}$ is odd.  We recall the
well-known fact that the 2-adic valuation of $\binom nm$ is equal to
the number of carries in the sum of $m$ and $n-m$ in base 2.  Since
$k$ is odd, exactly one of $k-1$ and $k-3$ is congruent to $2\pmod4$
and the other is congruent to $0\pmod4$.  To fix ideas, suppose $k-1$
is $0\pmod4$.  Let $a$ be an integer so that $2^a>k$ and choose $g$ so
that $2g-k-1$ is congruent to $0\pmod{2^a}$.  Then there are no
carries in the sum $(k-1)+(2g-k-1)$ and so $\binom{2g-2}{k-1}$ is odd.
On the other hand, both $k-3$ and $2g-k+1$ are congruent to $2\pmod4$
and so there is at least one carry in the sum $(k-3)+(2g-k+1)$ and
$\binom{2g-2}{k-3}$ is even.  The case where $k-1\equiv2\pmod4$ is
similar and will be left to the reader.  This completes the proof of
Theorem~\ref{thm:av2}.

\begin{rem}
  The above proof would carry over verbatim to $p=2$ if we had a curve
  over $\F_2(u)$ with large, everywhere tame monodromy, and inertia
  acting by unipotent reflections at an odd number of non-zero finite
  places.
\end{rem}

\section{Proof of Theorem~\ref{thm:ec1}}\label{s:ec1}
Let $p$ be any prime, $q$ a power of $p$, and $E$ an elliptic curve
over $F=\Fq(v)$.  Let $\rho$ be the representation of
$G_F$ on $H^1(E\times\overline{F},\Qlbar)$ for some
$\ell\neq p$; $\rho$ is self dual of weight 1 and sign $-1$.  For
every finite separable extension $K$ of $F$, we have
$L(E/K,s)=L(\rho,K,s)$.  We prove Theorem~\ref{thm:ec1} by showing
that after replacing $F$ with an extension of the form $\F_{r}(u)$,
the hypotheses of Theorem~\ref{thm:towers} are met and so $E$ obtains
large analytic rank over $\F_{r}(t)$ with $t^d=u$ for large $d$ of the
form $r^{n}+1$.

\begin{lemma}
  Suppose that $E$ is an elliptic curve over $F=\Fq(v)$ with
  $j(E)\not\in\Fq$.  Then there is a finite separable extension of $F$
  of the form $\Fr(u)$ over which $E$ has an $\Fr$-rational place of
  multiplicative reduction and two $\Fr$-rational places of good
  reduction.
\end{lemma}

\begin{proof}
  Extending the ground field to $\F_r$ and making a linear change of
  coordinates, we may assume that the $j$-invariant of $E$ has a pole
  at $v=0$.  It is well known that there is a finite separable
  extension of the completion $\F_r((v))$ over which $E$ obtains
  multiplicative reduction.

  The following lemma (whose proof is due to Bjorn Poonen) says that
  we may realize the local extension as the completion of a global
  extension of {\it rational\/} fields.  Admitting the lemma, $E$ has
  an $\Fr$-rational place of multiplicative reduction over $\F_r(u)$.
  Clearly at the expense of increasing $r$ we may insure that $E$ also
  has two $\Fr$-rational places of good reduction over $\F_r(u)$.
\end{proof}

\begin{lemma}
  Let $F=\F_r(v)$ and let $K_0$ be a finite separable extension of
  $F_0:=\F_r((v))$.  Then there exists a finite separable extension
  $K$ of $F$ of the form $K=\F_r(u)$ so that the completion of $K$ at
  the place $u=0$ is isomorphic, as extension of $F_0$, to $K_0$.
\end{lemma}

\begin{proof}
  It is well known that $K_0$ is abstractly isomorphic to
  $\F_r((\varpi))$.  Let $g(T)\in\F_r[[T]]$ ($T$ an indeterminate) be
  the formal series such that $g(\varpi)=v$ in $K_0$.  Since $K_0/F_0$
  is a separable extension, $g'(\varpi)$, the derivative series
  evaluated at $\varpi$, is not zero.  For a positive integer $n$, let
  $g_n$ be the sum of the first $n$ terms of $g$.  Then as
  $n\to\infty$, the valuation of $g_n'(\varpi)$ stabilizes at a finite
  value whereas the valuation of $g_n(\varpi)-v$ tends to infinity.
  By Hensel's lemma, for any sufficiently large $n$, there is a root
  $u$ of $g_n(T)-v$ in $K_0$ which is congruent to $\varpi$ modulo a
  high power of $\varpi$ and which is thus a uniformizer of $K_0$.
  Now let $K$ be the subfield of $K_0$ generated by $\F_r$ and $u$.
  Since $v$ is a polynomial in $u$ with a non-zero derivative, $K$ is
  a finite, separable extension of $F$.
\end{proof}

\subsection{}
With these preliminaries out of the way, we can prove
Theorem~\ref{thm:ec1}.  Since the $j$-invariant of $E$ is not in
$\Fq$, $\rho$ is irreducible.

If the degree of the conductor of $E$ is odd, we make a linear change
of coordinates so that $u=0$ and $u=\infty$ are places of good
reduction.  Then, in the notation of Theorem~\ref{thm:towers},
$\deg(\n')+\swan_0(\rho)+\swan_\infty(\rho)$ is the degree of the
conductor of $E$ which is odd.  Thus $\rho$ satisfies the hypotheses
of Theorem~\ref{thm:towers}.

If the degree of the conductor of $E$ is even, we make a change of
coordinates so that $u=0$ is a place of good reduction and $u=\infty$
is a place of multiplicative reduction.  Then, in the notation of
Theorem~\ref{thm:towers}, $\deg(\n')+\swan_0(\rho)+\swan_\infty(\rho)$
is one less than the degree of the conductor of $E$ and is therefore
odd.  Again $\rho$ satisfies the hypotheses of
Theorem~\ref{thm:towers}.  \qed

\subsection{}
In \cite{UlmerR3} we give examples of elliptic curves over $\Fq(u)$
with
bounded rank in the tower $\Fq(t)$ ($t^d=u$), $d\to\infty$.

\section{Proof of Theorem~\ref{thm:ec2}}\label{s:ec2}

We use the notation of the statement of Theorem~\ref{thm:ec2}.  Tate
and Shafarevitch observed in \cite{TS} that to produce a quadratic
twist $E'$ of $E$ with large rank, one must produce a hyperelliptic
curve $\Curve\to\P^1_t$ whose Jacobian has many factors isogenous to
$E_0$.  More precisely, if $E'$ is the twist of $E$ by the quadratic
extension $\Fp(\Curve)$, then the rank of $E'(\Fp(t))$ is equal to the
rank of the endomorphsim ring of $E_0$ (which in our case is 2) times
the number of isogeny factors of $J(C)$ isogenous to $E_0$.  (See
\cite{UlmerR3}*{\S4} for more details.)  Moreover, we may detect the
number of times a particular abelian variety appears in the Jacobian
of a curve via Honda-Tate theory by considering the inverse roots of
its zeta function.  In the rest of this section we will use an
orthogonal variant of the towers Theorem~\ref{thm:towers} to produce
hyperelliptic curves whose Jacobians have many isogeny factors
isogenous to a given supersingular elliptic curve.

\subsection{}
We call the inverse roots of the zeta function of a curve its {\it
  Weil numbers\/}.  It is well known (see, e.g.,
\cite{Waterhouse}*{Chap.~4}) that a supersingular elliptic curve over
$\Fp$ (any $p$) either has Weil number $\zeta_4\sqrt{p}$ with
$\zeta_4$ a primitive 4-th root of unity, or $p=3$ and the Weil number
is $\zeta_{12}\sqrt{3}$ with $\zeta_{12}$ a primitive 12-th root of
unity, or $p=2$ and the Weil number is $\zeta_{8}\sqrt{2}$ with
$\zeta_8$ a primitive 8-th root of unity.  We start with the case
$p>2$ and $E_0$ a supersingular elliptic curve with Weil number
$\zeta_4\sqrt{p}$.

\subsection{}\label{ss:FiberProducts}
Let $F=\Fp(u)$ and let $C_1$ be a geometrically irreducible curve
smooth and proper over $\Fp$ with genus $g\ge0$ equipped with a degree
2 morphism $\pi:C_1\to\P^1$.  We assume $\pi$ to be ramified at $2g+2$
geometric points ($a_1,\dots,a_{2g+1},\infty$) one of which is
infinity and none of which are 0.  Corresponding to the covering
$\pi:C_1\to\P^1$ is a character $\rho$ of $G_F$ of
order 2.  The numerator of the zeta function of $C_1$ is
$L(\rho,F,s)$.

For a positive integer $d$ not divisible by $p$, we let $F_d=\Fp(t)$
be the extension of $F$ with $u=t^d$ and we let $\pi_d:C_d\to\P^1_t$
be the 2-1 covering corresponding to $\rho$ restricted to
$\gal(\overline{F}/\Fp(t))$.  It is not hard to check that $C_d$ is
the normalization of the fiber product $C_1\times_{\P^1_u}\P^1_t$.
The ramification points of $\pi_d$ are the $d$-th roots of the $a_i$,
0, and, if $d$ is odd, $\infty$.

\subsection{}
The numerator of the zeta function of $C_d$ is $L(\rho,F_d,s)$ and by
the analysis in Subsection~\ref{ss:geom-abelian-base-change},
\begin{equation*}
L(\rho,F_d,s)=\prod_{o\subset\Z/d\Z}L(\rho\tensor\sigma_o,F,s)
\end{equation*}
where the product is over the orbits of multiplication by $p$ on
$\Z/d\Z$.  If $\chi$ is a character of
$\gal(\overline{F}/\Fp(\mu_d,u))$ of order $d$ corresponding to the extension
$\Fp(\mu_d,t)$, then it is easy to see that
$\deg\cond(\rho\tensor\chi^i)$ is odd except when $\chi^i$ has order
dividing 2.  It follows from Theorem~\ref{thm:Lzeroes} that if
$d=p^n+1$ and $o\subset\Z/d\Z$ is any orbit for multiplication by $p$
other than $\{0\}$ or $\{d/2\}$ and $a=\#o$, then $1+(Tp^{1/2})^a$
divides $L(\rho\tensor\sigma_o,F,s)$.  If we take $n$ odd and let $o$
be an orbit passing through $i\in(\Z/d\Z)^\times$, then $a$ is $2n$ and so
$\zeta_4\sqrt{p}$ is a root of $1+(Tp^{1/2})^a$.  Since there are
$\phi(p^n+1)/2n$ such orbits, we see that $\zeta_4\sqrt{p}$ is an
inverse root of the numerator of the zeta function of $C_d$ with large
multiplicity.  As discussed above, this shows that $E'$ has large rank
over $F_d=\Fp(t)$.

\subsection{}
Now consider the case where $E_0$ is a supersingular elliptic curve
over $\F_3$ with Weil number $\zeta_{12}\sqrt{3}$ where $\zeta_{12}$
is a primitive 12-th root of unity.  We proceed as above except that
we assume that $n\equiv 3\pmod{6}$ so that $\zeta_{12}\sqrt{3}$ is an
inverse root of $1+(T\sqrt{3})^a$ where $a=2n$.

\subsection{}
If $p=2$ then we proceed as above starting with a curve $C_1\to\P^1_u$
corresponding to a quadratic character $\rho$ satisfying the conductor
condition of \ref{thm:towers}, namely that
$\swan_0(\rho)+\swan_\infty(\rho)+\deg\n'$ is odd.  (For example,
$y^2+y=u$.)  If the Weil number of $E_0$ is $\zeta_4\sqrt{2}$ then we
take $d=p^n+1$ with $n$ odd and if the Weil number of $E_0$ is
$\zeta_8\sqrt{2}$ then we take $d=p^n+1$ with $n\equiv2\pmod4$.

\subsection{}
Interestingly, the argument above does not prove that if $E_0$ is any
supersingular elliptic curve over $\Fq$ then there are quadratic
twists of $E$ with high rank, only the slightly weaker statement that
there is a power $r$ of $q$ and quadratic twists of $E$ with high rank
over $\F_r(t)$.  The problem is that if the Weil number of $E_0$ is
$\zeta_m\sqrt{q}$ with $m$ odd, then this Weil number is not a root of
$1+(q^{1/2}T)^a$ for any even $a$.

\numberwithin{equation}{section}
\section{A remark on rank bounds}\label{s:RankBounds}
Suppose as usual that $F$ is the function field of a curve $\Curve$ of
genus $g_\Curve$ over $\Fq$ and that $\rho$ is a representation of
$G_F$ satisfying the hypotheses of
Subsection~\ref{ss:rho-hyps} and (for simplicity) that $\rho$ restricted to
$\gal(\overline{F}/\Fqbar F)$ has no trivial constituents.  Let $\n$ be
the conductor of $\rho$.  Then the Grothendieck-Ogg-Shafarevitch
formula says that the degree of the $L$-function $\rho$ over $F$ as a
polynomial in $q^{-s}$ is $D=\deg(\n)+\deg(\rho)(2g_\Curve-2)$.  In
particular, we have the ``geometric'' rank bound (cf.~\cite{UlmerA})
\begin{equation}\label{eqn:geometric-bound}
\ord_{s=(w+1)/2}L(\rho,F,s)\le\ord_{s=(w+1)/2}L(\rho,\F_r F,s)\le D
\end{equation}
valid for any power $r$ of $q$.

This can be improved when $D$ is large with respect to $q$,
$g_\Curve$, and $\deg(\rho)$.  Indeed, minor modifications of Brumer's
argument in \cite{Brumer} (itself modelled on Mestre's \cite{Mestre})
allow one to prove the arithmetic rank bound
\begin{equation}\label{eqn:arithmetic-bound}
\ord_{s=(w+1)/2}L(\rho,F,s)\le\frac{D}{2\log_q D}
+O\left(\frac{D}{(\log_qD)^2}\right)
\end{equation}
where the implied constant depends only on $q$, 
$g_\Curve$, and $\deg(\rho)$.

In \cite{UlmerR} we showed that the main term of this arithmetic
bound, as well as the geometric bound, are sharp for $L$-functions of
elliptic curves.  The towers Theorem~\ref{thm:towers} gives a large
supply of other examples related to this question.

Indeed, suppose that $\rho$ is a representation of
$G_F$ where $F=\Fq(u)$ satisfying the hypotheses of
\ref{thm:towers} and let $N$ be the quantity
$\swan_0(\rho)+\swan_\infty+\deg(\n')$ appearing in that result.  Then
the degree of $L(\rho,F_d,s)$, where $F_d=\Fq(t)$ with $u=t^d$, is
asymptotic to $Nd$; they differ by an amount bounded independently of
$d$.  The towers Theorem~\ref{thm:towers} shows that for $d$ of the
form $d=q^n+1$
$$\ord_{s=(w+1)/2}L(\rho,F_d,s)\ge\frac{d}{2\log_q d}-c$$
and
$$\ord_{s=(w+1)/2}L(\rho,\F_{q^{2n}}F_d,s)\ge d-c'$$
where $c$ and $c'$ are constants independent of $n$.  These lower
bounds are roughly $1/N$ times the upper bounds discussed above.

For the curves discussed in Section~\ref{s:algranks} with $p>2$, we
have $N=1$ and so we have a large collection of interesting
representations for which the main term of the rank bounds are sharp.


\end{document}